\documentclass{gtart_h}

\def\ifplaintex{\expandafter\ifx\csname documentclass\endcsname\relax}

\def\gtp{{\mathsurround=0pt\it $\cal G\mskip-2mu$eometry \&\ 
$\cal T\!\!$opology $\cal P\!$ublications}}  

\def\recd{{\small Received:\qua\receiveddate\ifx\reviseddate\relax
\else\qquad Revised:\qua\reviseddate\fi\par}} 

\def\lognumber#1{\def\thelognumber{#1}}
\def\volumenumber#1{\def\thevolumenumber{#1}}
\def\volumeyear#1{\def\thevolumeyear{#1}}
\def\papernumber#1{\def\thepapernumber{#1}}
\def\pagenumbers#1#2{\def\startpage{#1}\def\finishpage{#2}}
\def\published#1{\def\publishdate{#1}}

\def\received#1{\def\receiveddate{#1}}
\def\revised#1{\def\reviseddate{#1}}
\def\accepted#1{\def\accepteddate{#1}}
\def\asciititle#1{\def\theasciititle{#1}}

\def\asciiauthors#1{\def\theasciiauthors{#1}}
\def\asciiaddress#1{\def\theasciiaddress{#1}}
\def\asciiemail#1{\def\theasciiemail{#1}}

\def\coverauthors#1{\def\thecoverauthors{#1}}
\long\def\asciiabstract#1{\long\def\theasciiabstract{#1}}

\let\\\par\let\thelognumber\relax\let\thevolumenumber\relax
\let\thepapernumber\relax\let\thevolumeyear\relax\let\startpage\relax
\let\finishpage\relax\let\publishdate\relax\let\receiveddate\relax
\let\reviseddate\relax\let\accepteddate\relax\let\theasciititle\relax
\let\theasciiauthors\relax\let\theasciiaddress\relax
\let\theasciiabstract\relax

\let\thecoverauthors\relax\let\theasciiemail\relax

\ifplaintex
\font\logobig=cmssbx10 scaled 3836
\font\logomed=cmssbx10 scaled 2557
\else
\font\logobig=cmssbx10 scaled 4200
\font\logomed=cmssbx10 scaled 2800
\fi

\long\def\makeagttitle{   
\count0=\startpage
\agt\hfill      
\hbox to 45truept{\vbox to 0pt{\vglue -13truept{\logomed A\kern -.37em{\logobig 
T}\kern -.38em G}\vss}\hss}
\break
{\small Volume \thevolumenumber\ (\thevolumeyear)
\startpage--\finishpage\nl
Published: \publishdate}

\vglue .25truein

{\parskip=0pt\leftskip 0pt plus
1fil\def\\{\par\smallskip}{\Large\bf\thetitle}\par\medskip} \vglue
0.05truein

{\parskip=0pt\leftskip 0pt plus 1fil\def\\{\par}{\sc\theauthors}
\par\medskip}%
 
\vglue 0.03truein 

{\small\leftskip 25truept\rightskip 25truept{\bf Abstract}\stdspace\theabstract

{\bf AMS Classification}\stdspace\theprimaryclass
\ifx\thesecondaryclass\relax\else; \thesecondaryclass\fi\par
{\bf Keywords}\stdspace \thekeywords\par}\vglue 7truept

}   

\ifplaintex
\font\phead=cmsl9 scaled 950
\font\pnum=cmbx10 scaled 913
\font\pfoot=cmsl9 scaled 950
\headline{\vbox to 0pt{\vskip -4.5mm\line{\small\phead\ifnum
\count0=\startpage ISSN 1472-2739 (on-line) 1472-2747 (printed)
\hfill {\pnum\folio}\else\ifodd\count0\def\\{ }%
\ifx\theshorttitle\relax\thetitle\else\theshorttitle\fi\hfill{\pnum\folio}
\else\def\\{ and }{\pnum\folio}\hfill\ifx\theshortauthors\relax\theauthors
\else\theshortauthors\fi\fi\fi}\vss}}
\footline{\vbox to 0pt{\vglue 0mm\line{\small\pfoot\ifnum\count0=\startpage
\copyright\ \gtp\hfill\else
\agt, Volume \thevolumenumber\ (\thevolumeyear)\hfill\fi}\vss}}
\else
\font=cmsl9 scaled 1050
\font\lnum=cmbx10 
\font=cmsl9 scaled 1050
\makeatletter
\def\@oddhead{{\small\lhead\ifnum\count0=\startpage ISSN 1472-2739 
(on-line) 1472-2747 (printed)\hfill {\lnum\number\count0}\else\ifodd\count0
\def\\{ }\ifx\theshorttitle\relax \thetitle \else\theshorttitle\fi\hfill
{\lnum\number\count0}\else\def\\{ and }{\lnum\number\count0}
\hfill\ifx\theshortauthors\relax 
\theauthors\else\theshortauthors\fi\fi\fi}}\def\@evenhead{\@oddhead}
\def\@oddfoot{\small\lfoot\ifnum\count0=\startpage\copyright\ \gtp\hfill\else
\agt, Volume \thevolumenumber\ (\thevolumeyear)\hfill\fi}
\def\@evenfoot{\@oddfoot}
\makeatother
\fi
\let\maketitlepage\makeagttitle
\let\makeshorttitle\maketitlepage
\let\maketitle\maketitlepage

\newwrite\gtoutfile
\long\gdef\makeheadfile{  
{\def\\{, }\def\s{ }
\immediate\openout\gtoutfile head.xxx
\immediate\write\gtoutfile{Proxy-for: \ifx\theasciiauthors\relax
\theauthors\else\theasciiauthors\fi\s<\ifx\theasciiemail\relax\theemail\else\theasciiemail\fi>}
\immediate\write\gtoutfile{\noexpand\\}
\immediate\write\gtoutfile{Authors: \ifx\theasciiauthors\relax
\theauthors\else\theasciiauthors\fi}
{\def\\{ }\immediate\write\gtoutfile{Title: \ifx\theasciititle\relax
\thetitle\else\theasciititle\fi}}
\immediate\write\gtoutfile{Subj-class: GT or SG, GR etc}
\immediate\write\gtoutfile{MSC-class: \theprimaryclass\ifx\thesecondaryclass\relax\else, \thesecondaryclass\fi}
\immediate\write\gtoutfile{Journal-ref: Algebr. Geom. Topol. \thevolumenumber\s
(\thevolumeyear) \startpage-\finishpage}
\immediate\write\gtoutfile{Comments: Published by Algebraic and
Geometric Topology at}
\immediate\write\gtoutfile{\s\s\s  http://www.maths.warwick.ac.uk/agt/AGTVol\thevolumenumber/agt-\thevolumenumber-\thepapernumber.abs.html}
\immediate\write\gtoutfile{\noexpand\\}
\immediate\write\gtoutfile{}
\ifx\theasciiabstract\relax
\immediate\write\gtoutfile{\theabstract}\else
\immediate\write\gtoutfile{\theasciiabstract}\fi
\immediate\write\gtoutfile{}
\immediate\write\gtoutfile{\noexpand\\}
\immediate\write\gtoutfile{}
\immediate\closeout\gtoutfile}}  

\def\maketitlepage{\makeagttitle\makeheadfile}
\let\makeshorttitle\maketitlepage
\let\maketitle\maketitlepage

\lognumber{52}
\volumenumber{4}
\volumeyear{2004}
\papernumber{52}
\pagenumbers{1177}{1210}
\received{23 September 2004} 
\revised{6 December 2004}
\accepted{6 December 2004}
\published{15 December 2004}

\usepackage{epsfig}

\usepackage{amssymb,amsmath}

\newtheorem{theorem}{Theorem}[section]
\newtheorem{lemma}[theorem]{Lemma}
\newtheorem{corollary}[theorem]{Corollary}
\newtheorem{proposition}[theorem]{Proposition}

\newtheorem{example}[theorem]{Example}

\newcommand{\al}{\alpha}
\newcommand{\be}{\beta}

\renewcommand{\dh}{\hat d}
\newcommand{\RP}{{\mathbb R}{\mathbb P}}
\newcommand{\R}{\mathbb R}
\newcommand{\C}{\mathcal C}
\renewcommand{\S}{\mathcal S}
\newcommand{\Z}{\mathbb Z}
\newcommand{\B}{\mathcal B}

\newcommand{\N}{\mathbb N}
\newcommand{\lp}{{\Z[A^{\pm 1}]}}

\newcommand{\comment}[1]{}
\newcommand{\pic}[1]{\parbox{.6cm}{\psfig{figure=#1.eps,height=.6cm}}}

\newcommand{\picll}[1]{\parbox{.8cm}{\psfig{figure=#1.eps,height=.8cm}}}

\def\down#1{\raise-3pt\hbox{#1}}

\newcommand{\cp}{\down{\begin{picture}(12,12) \put(6,6){\circle{12,12}}
\put(6,6){\makebox(0,0){$\scriptstyle +$}} \end{picture}\hspace*{1pt}}}

\newcommand{\cm}{\down{\begin{picture}(12,12) \put(6,6){\circle{12,12}}
\put(6,6){\makebox(0,0){-}} \end{picture}\hspace*{1pt}}}

\newcommand{\ce}{\down{\begin{picture}(12,12) \put(6,6){\circle{12,12}}
\put(6,6){\makebox(0,0){$\scriptstyle 
\varepsilon$}}\end{picture}\hspace*{1pt}}}

\newcommand{\cme}{\down{\begin{picture}(12,12)\put(6,6){\oval(12,12)}
\put(6,6){\makebox(0,0){-$\scriptstyle \varepsilon$}}\end{picture}
\hspace*{.2pt}}}

\newcommand{\cpe}{\down{\begin{picture}(12,12)\put(6,6){\oval(12,12)}
\put(6,6){\makebox(0,0){+$\scriptstyle \varepsilon$}}\end{picture}
\hspace*{.2pt}}}

\newcommand{\cpo}{\down{\begin{picture}(12,12) \put(6,6){\oval(12,12)}
\put(6,6){\makebox(0,0){$\scriptstyle +0$}} \end{picture}\hspace*{1pt}}}

\newcommand{\cmo}{\down{\begin{picture}(12,12)\put(6,6){\oval(12,12)}
\put(6,6){\makebox(0,0){-$\scriptstyle 0$}} \end{picture}\hspace*{1pt}}}

\newcommand{\ceo}{\down{\begin{picture}(12,12)\put(6,6){\oval(12,12)}
\put(6,6){\makebox(0,0){$\scriptstyle \varepsilon 0$}}\end{picture}}}

\newcommand{\three}[3]{\renewcommand{\arraystretch}{0.2}
\begin{array}{ccc} #1\\ #2\\ #3 \end{array}}

\newcommand{\diagtwoa}[8]{\renewcommand{\arraystretch}{0.2}
\begin{array}{c@{\,}c@{}c@{}c@{}c}
 & &{#2\atop #3}& &\\
 &\nearrow & & \searrow & \\
\scriptstyle #1& & & & \begin{array}{c} \scriptstyle #4\\ \scriptstyle
 #5\\ \scriptstyle #6\end{array}\\
 &\searrow & & \nearrow & \\
 & &{#7\atop #8}& &
\end{array}}

\newcommand{\diagtwob}[8]{\renewcommand{\arraystretch}{0.2}
\begin{array}{c@{\,}c@{}c@{}c@{}c}
 & &{#2\atop #3}& &\\
 &\swarrow & & \searrow & \\
\scriptstyle #1& & & & \begin{array}{c} \scriptstyle #4\\ \scriptstyle
 #5\\ \scriptstyle #6\end{array}\\
 &\searrow & & \swarrow & \\
 & &{#7\atop #8}& &
\end{array}}

\newcommand{\diagfoura}[6]{\renewcommand{\arraystretch}{0.2}
\begin{array}{c@{\,}c@{}c@{}c@{}c}
 & &{#2\atop #3}& & \\
 &\nearrow & & \searrow & \\
\scriptstyle #1& & & & {#4\atop #5}\\
 &\searrow & & \nearrow & \\
 & & \scriptstyle #6 & &
\end{array}}

\renewcommand{\t}{$\tau$}
\newcommand{\ve}{\varepsilon}

\begin{document}

\title[Categorification of the Kauffman bracket skein
  module]{Categorification of the Kauffman bracket skein module of 
$I$-bundles over surfaces}

\asciititle{Categorification of the Kauffman bracket skein module of 
I-bundles over surfaces}

\authors{Marta M. Asaeda\\J\'ozef H. Przytycki\\Adam S. Sikora}
\asciiauthors{Marta M. Asaeda\\Jozef H. Przytycki\\Adam S. Sikora}
\coverauthors{Marta M. Asaeda\\J\noexpand\'ozef H. Przytycki\\Adam S. Sikora}
\shortauthors{Asaeda, Przytycki and Sikora}

\address{Dept of Mathematics, 14 MacLean Hall\\University of Iowa,  
Iowa City, IA 52242, USA\\\medskip\\
Dept of Mathematics, Old Main Bldg, The George 
Washington University\\1922 F St NW, Washington, DC 20052, USA\\
\medskip\\
Dept of Mathematics, 244 Math Bldg, SUNY at Buffalo\\ 
Buffalo, NY 14260, USA\quad {\rm and}\\
Inst for Adv Study, School of Math, Princeton, NJ 08540, USA}
\asciiaddress{Dept of Mathematics, 14 MacLean Hall\\University of Iowa,  
Iowa City, IA 52242, USA\\and\\
Dept of Mathematics, Old Main Bldg, The George 
Washington University\\1922 F St NW, Washington, DC 20052, USA\\
and\\
Dept of Mathematics, 244 Math Bldg, SUNY at Buffalo\\ 
Buffalo, NY 14260, USA, and\\
Inst for Adv Study, School of Math, Princeton, NJ 08540, USA}
\gtemail{\mailto{asaeda@math.uiowa.edu}, 
\mailto{przytyck@gwu.edu}, \mailto{asikora@buffalo.edu}}
\asciiemail{asaeda@math.uiowa.edu, przytyck@gwu.edu, asikora@buffalo.edu}

\begin{abstract} 
Khovanov defined graded homology groups for links $L\subset \R^3$
and showed that their polynomial Euler characteristic
is the Jones polynomial of $L$.
Khovanov's construction does not extend in a straightforward
way to links in $I$-bundles $M$ over surfaces $F\ne D^2$
(except for the homology with $\Z/2$ 
coefficients only). Hence, the goal of this paper is to provide
a nontrivial generalization of his method leading to homology 
invariants of links in $M$ with arbitrary rings of coefficients.

After proving the invariance of our homology groups under Reidemeister
moves, we show that the polynomial Euler characteristics of our
homology groups of $L$ determine the coefficients of $L$ in the
standard basis of the skein module of $M.$ Therefore, our homology groups
provide a ``categorification'' of the Kauffman bracket skein module of $M.$
Additionally, we prove a generalization of Viro's exact sequence 
for our homology groups.
Finally, we show a duality theorem relating cohomology 
groups of any link $L$ to the homology groups of the mirror image of $L.$
\end{abstract}

\asciiabstract{%
Khovanov defined graded homology groups for links L in R^3 and showed
that their polynomial Euler characteristic is the Jones polynomial of
L.  Khovanov's construction does not extend in a straightforward way
to links in I-bundles M over surfaces F not D^2 (except for the
homology with Z/2 coefficients only).  Hence, the goal of this paper
is to provide a nontrivial generalization of his method leading to
homology invariants of links in M with arbitrary rings of
coefficients.  After proving the invariance of our homology groups
under Reidemeister moves, we show that the polynomial Euler
characteristics of our homology groups of L determine the coefficients
of L in the standard basis of the skein module of M.  Therefore, our
homology groups provide a `categorification' of the Kauffman bracket
skein module of M.  Additionally, we prove a generalization of Viro's
exact sequence for our homology groups.  Finally, we show a duality
theorem relating cohomology groups of any link L to the homology
groups of the mirror image of L.}



\maketitle

\section{Introduction}

In his seminal work \cite{K1}, Khovanov constructed
graded homology groups for links $L\subset \R^3$ and proved that
their polynomial Euler characteristic is the 
Jones polynomial of $L.$ 
(Throughout this paper we use Viro's (framed) normalization of
Khovanov homology groups.)
The crucial element in Khovanov's construction is the representation
of a link by its diagram in $D^2$ unique up to 
Reidemeister moves. Despite the fact that links in all
$I$-bundles over surfaces allow analogous diagrammatic representations, 
Khovanov's construction does not extend in a straightforward way to
such links. In this paper we overcome this difficulty
and define homology invariants of links in all orientable $I$-bundles $M$
over surfaces $F\ne \RP^2.$

We say that a simple closed loop $\gamma\subset F$ is {\em bounding} 
if it bounds either a disk or a M\"obius band in $F.$
Let ${\C}(F)$ be the set of all unoriented, unbounding, 
simple closed curves in $F$ considered up to homotopy. 
In Sections \ref{s_chain}-\ref{s_d^2}, we construct 
homology groups, $H_{ijs}(L),$ of links $L\subset M$ for any 
$i,j\in \Z,$ and $s\in \Z\C(F),$
and we prove (Section \ref{main_statements}) that these
homology groups are invariant under isotopies of $L$ in $M.$

In Section \ref{s_viro}, we construct a generalization
of Viro's exact sequence, \cite[Section 3.3]{V}, which associates with
any skein triple of link diagrams in $F$ a long exact sequence of 
their homology groups. 

In Section \ref{s_categorify}, we prove that if $F$ is orientable then 
the polynomial Euler characteristics of our homology groups
$H_{ijs}(L)$ determine the coefficients of $L$ in the standard basis of 
the Kauffman bracket skein module. Additionally, we
prove an analogous result for links in twisted $I$-bundles over
unorientable surfaces.

In Section \ref{cohomology}, we introduce cohomology groups of links
and prove a duality theorem relating them to homology groups 
of the mirror image of $L.$  Finally, we discuss 
the dependence of our homology groups of links on the $I$-bundle structure 
of the ambient $3$-manifold in Section \ref{I-bundle}.

The methods of this paper are applied in \cite{APS} to define a new 
homology theory for tangles. 

{\bf Acknowledgments}\qua
The first author was sponsored in part by NSF grant
DMS-0202613. The second author was partially sponsored by grant 
TMP-25051958. The third author was sponsored in part by NSF grants
DMS-0307078 and DMS-0111298.

%
\section{Skein modules of band links}
\label{s_skein_mod}
%

Let $M$ be an orientable $I$-bundle over a surface $F.$
Consequently,
$M$ is either $F\times I$ for an orientable $F,$ or $M$ is the twisted 
$I$-bundle over an unorientable surface $F.$ 
In the latter case we assume that $F\ne \RP^2.$ 
A {\em band knot} $K$ in $M$ is either
\begin{itemize}
\item an embedding of an annulus into $M$
such that the projection of its core into $F$ is a preserving
orientation loop in $F$; or
\item an embedding of a M\"obius band into $M$
such that the projection of its core into $F$ is a reversing orientation 
curve in $F.$
\end{itemize}
A disjoint union of band knots is called a {\em band link}.
Each band link in $M$ is represented by a diagram in $F.$
(All link diagrams are considered with their blackboard framing.
Hence, an orientation reversing loop in $F$ represents a M\"obius
band in $M.$) From now on, all links considered in this paper are band 
links, unless stated otherwise. 

Denote the set of all band links in an orientable $I$-bundle $M$ over
$F,$ by ${\mathcal L}_b(M).$ For a given ring $R$ with a distinguished
element $A^{\pm 1},$ the skein module of band links in $M$, 
$\S_b(M;R),$ is the quotient
of $R{\mathcal L}_b(M)$ by the standard Kauffman bracket skein
relations: $$\pic{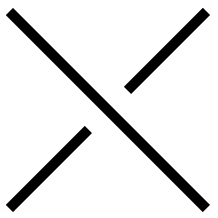}=A\pic{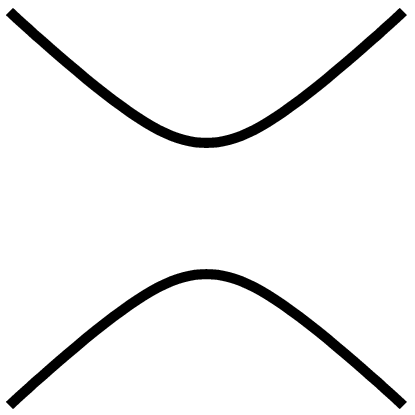}+A^{-1}\pic{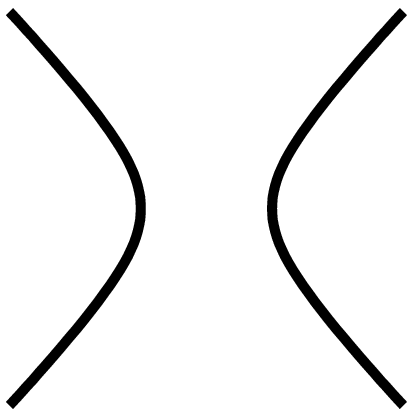},\quad 
L\cup \bigcirc =-(A^2+A^{-2})L.$$
Although the definition of a band link $L\subset M$ depends on the $I$-bundle
structure of $M,$ the theorem below shows the skein module $\S_b(M;R)$
(considered up to an isomorphism of $R$-modules) does not.
Let $\B(F)$ (respectively: $\B_{nb}(F)$) be the set of all link
diagrams in $F$ with no crossings and with no trivial (respectively: 
no bounding)
components. Both $\B(F)$ and $\B_{nb}(F)$ contain the empty link, $\emptyset.$
Following the proof of \cite[Theorem 3.1]{Pr}, one shows:

\begin{theorem}\label{skeinmod}
$\S_b(M;R)$ is a free $R$-module with a basis composed by (band) links
  represented by diagrams in $\B(F).$
\end{theorem}

Consequently, $\S_b(M;\Z[A^{\pm 1}])=\S(M;\Z[A^{\pm 1}])$ as
$R$-modules despite the fact that there is no obvious explicit
isomorphism between these modules for unorientable $F.$
On the other hand, there is a natural isomorphism between 
$\S_b(M;R)$ and $\S(M;R)$ for any ring $R$ containing $\sqrt{-A}.$
For such $R,$ we have an isomorphism $\lambda: 
\S_b(M;\Z[A^{\pm 1}])\to \S(M;\Z[A^{\pm 1}])$ sending 
$L=K_1\cup ...\cup K_n$ to $(-A)^{3k(L)/2}K_1'\cup ...\cup K_n',$ where
$K_i'=K_i$ if $K_i$ is an annulus and, otherwise, $K_i'$ is obtained
from $K_i$ by
adding a negative half-twist to $K_i$. Here $k(L)$ denotes the
number of M\"obius bands among the components, $K_i,$ of $L.$

%
\section{Chain groups}
\label{s_chain}
%

Let $D$ be a link diagram in $F$.
Following \cite{V}, a {\em Kauffman state}
of $D$ is an assignment of $+1$ or $-1$ marker to each of the
crossings of $D.$\footnote{ J.B.Listing 
in 1847 \cite{Li} was the first 
to consider two types of markers of a crossing: dexiotropic and
laeotropic.} An {\em enhanced (Kauffman) state} of 
$D$ is a Kauffman state with an additional assignment of 
$+$ or $-$ sign to each closed loop obtained 
by smoothing the crossings of $D$
according to the following convention:\vspace*{.1in}

\centerline{\psfig{figure=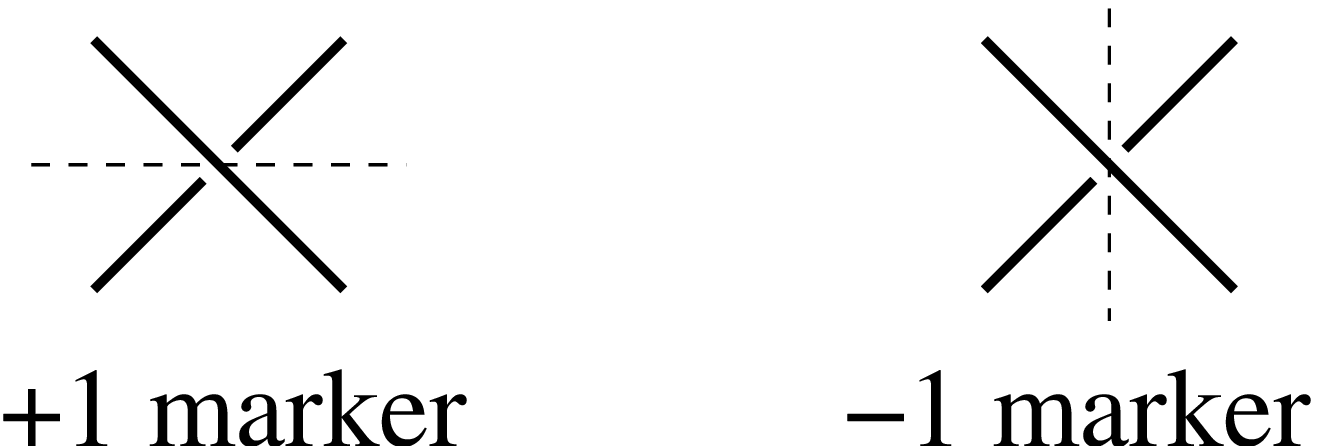,height=1.5cm}}

Such closed loops will be called by us ``circles''.

Denote the set of enhanced states of $D$ by $\S(D).$
The following notation is used throughout the paper for any $S\in \S(D):$ 
\begin{eqnarray*}
I(S) &=& \sharp \{\text{positive markers}\} - \sharp \{ 
\text{negative markers}\}.\\
J(S) &=&I(S)+2 \tau(S), \text{where}\\
\tau(S) &=& \sharp \{\text{positive trivial circles} \} -\sharp \{
\text{negative trivial circles}\}.
\end{eqnarray*}
A circle is trivial if it bounds a disk in $F$.
Furthermore, we say that a circle is bounding if it bounds either a
disk or a M\"obius band in $F.$
Let ${\C}(F)$ be the set of all unoriented, unbounding, 
simple closed curves in $F$ considered up to homotopy. 
If the unbounding components of an enhanced state $S$ 
are $\gamma_1,...,\gamma_n$ (some of which may be parallel to each other, 
and hence equal in $\C(F)$) and if these closed curves are
marked by $\ve_1,...,\ve_n\in \{+1,-1\},$ then
let
$$\Psi(S)=\sum_i \ve_i \gamma_i \in \Z{\C}(F).$$
Let $\S_{ijs}(D)$ be the set of enhanced states $S$ of $D$ with 
$i(S)=i,j(S)=j,$ and $\Psi(S)=s,$ and let $C_{ijs}(D)$ be the free 
abelian group spanned by enhanced states in $\S_{ijs}(D).$ 
Note that $C_{**s}(D)=0$ if $s=\sum_{i=1}^n \ve_i \gamma_i$ and
the curves $\gamma_i,$ $i=1,...,n,$ cannot be placed disjointly in
$F.$ On the other hand, if the curves $\gamma_i,$ $i=1,...,n,$ can be 
placed disjointly in $F,$ then $s$ represents a unique element of 
$\B(F).$ Therefore, alternatively, the index $s\in \Z{\C}(F)$ can be
replaced by $b\in \B(F).$

%
\section{Differentials}
\label{s_diff}
%

The definitions of the Khovanov chain groups of a
link in $\R^3$ given in \cite{K1,BN1,V} depend only on the choice of its
diagram $D.$ 
However, the definitions of differentials require
an ordering of crossings of $D.$
The same feature appears in our construction. Therefore, assume that 
the crossings of a link diagram $D$ in $F$ are ordered.
Let the differential
$$d_{ijs}: C_{ijs}(D) \longrightarrow C_{i-2,j,s}(D),$$
be defined for enhanced states $S\in C_{ijs}(D)$ as follows
$$d_{ijs}(S)= \sum_{v} (-1)^{t(S,v)} d_v(S), $$
where the sum is over all crossings $v$ in $D$ and the 
``partial derivative" in the ``direction" of $v$ is 
$$d_v(S) = \sum_{S'\in {\S}(D)} 
\hspace*{-.1in} [S:S']_v\, S'.$$ The incidence number, $[S:S']_v,$ 
is $1$ if the following four conditions are satisfied:

\begin{itemize}
\item[(a)] the crossing $v$ is marked by $+$ in $S$ and by $-$ in $S'$,

\item[(b)] $S$ and $S'$ assign the same markers to all the other crossings,

\item[(c)] the labels of the common circles in $S$ and $S'$ are unchanged,

\item[(d)] $J(S)=J(S'),$ $\Psi(S)=\Psi(S').$
\end{itemize}

Otherwise $[S:S']_v$ is equal to $0$.

The symbol $t(S,v)$ denotes the number of negative markers
assigned to crossings in $S$ bigger than $v.$ 

For fixed $D,$ $i,j,s,$ we write
$$\dh_{v}(S) =(-1)^{t(S,v)} d_v(S).$$
Hence $$d_{ijs}(S)=\sum_{{\rm crossings \ } v {\rm \ in} \ D} \hspace*{-.15in}
\dh_{v}(S)$$ for any $S \in {\S}_{ijs}(D).$ 

If $[S:S']_v=1$ then $\tau(S) = \tau(S')+1,$ 
and, therefore, $S$ and $S'$ in a neighborhood of the crossing 
$v$ have one of the following forms, where $\varepsilon=+$~or~$-$. 
\begin{center}
$\begin{array}{||ccc||}
\hline \vrule width 0pt depth 5pt height 12pt S & \Rightarrow & S' \\
\hline {\psfig{figure=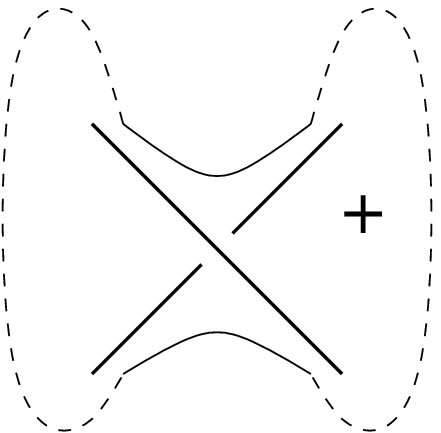,height=1.5cm}}
& &
{\psfig{figure=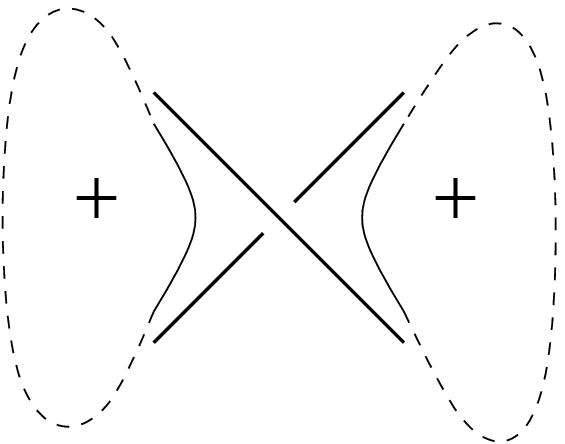,height=1.5cm}} 
\\
\hline 
{\psfig{figure=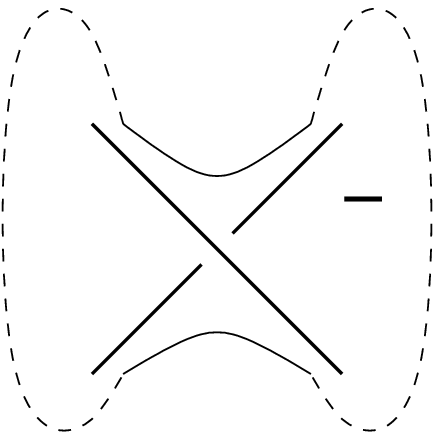,height=1.5cm}}
& &
{\psfig{figure=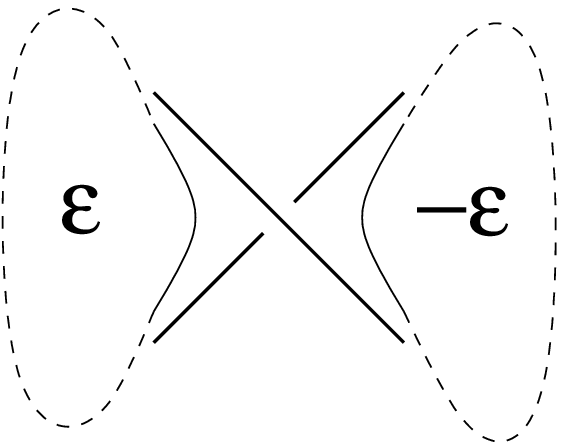,height=1.5cm}}
\\
\hline 
{\psfig{figure=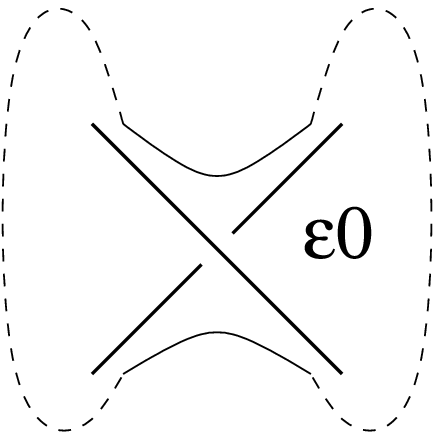,height=1.5cm}}
&  &
{\psfig{figure=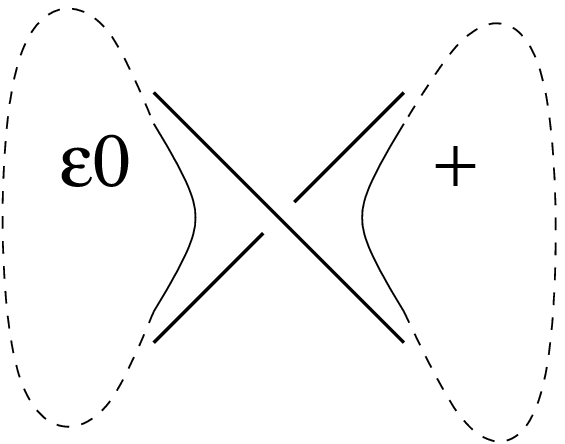,height=1.5cm}}
\\ 
\hline {\psfig{figure=ponem.eps,height=1.5cm}}& &
{\psfig{figure=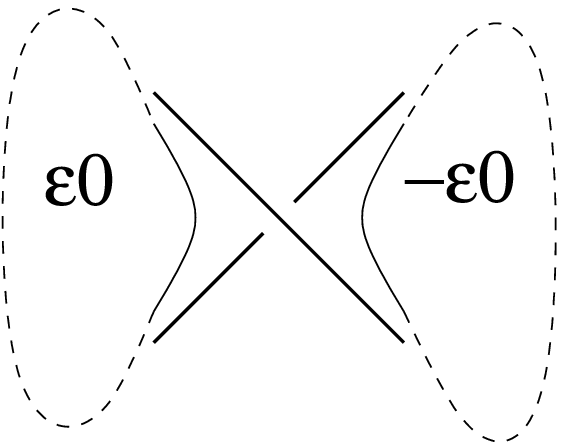,height=1.5cm}}\\ 
\hline 
{\psfig{figure=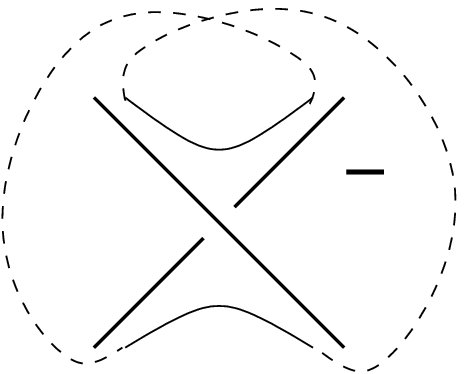,height=1.5cm}} 
& & 
{\psfig{figure=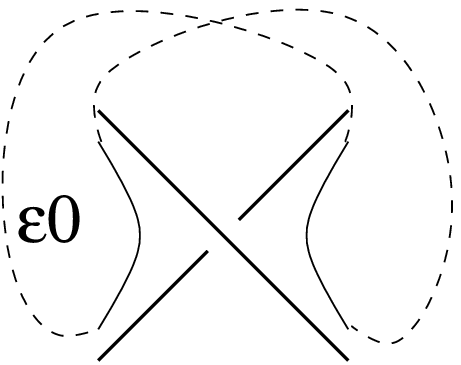,height=1.5cm}}
\\ 
\hline
\end{array}$
$\begin{array}{||ccc||}
\hline S\vrule width 0pt depth 5pt height 12pt & \Rightarrow & S'\\ \hline
{\psfig{figure=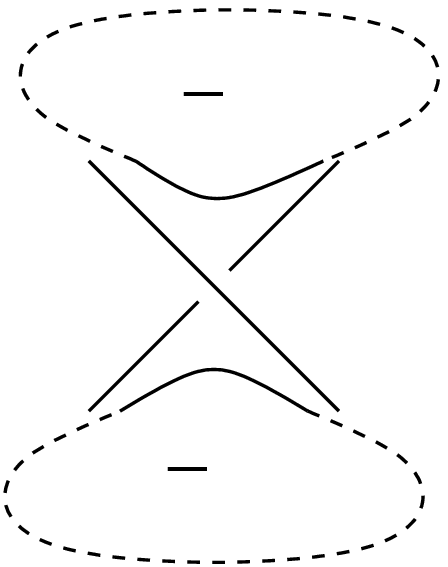,height=1.5cm}}& &
{\psfig{figure=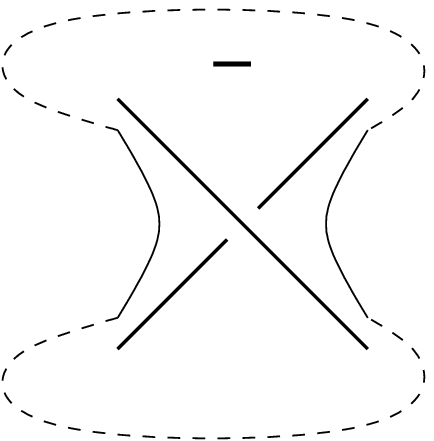,height=1.5cm}}\\ \hline
{\psfig{figure=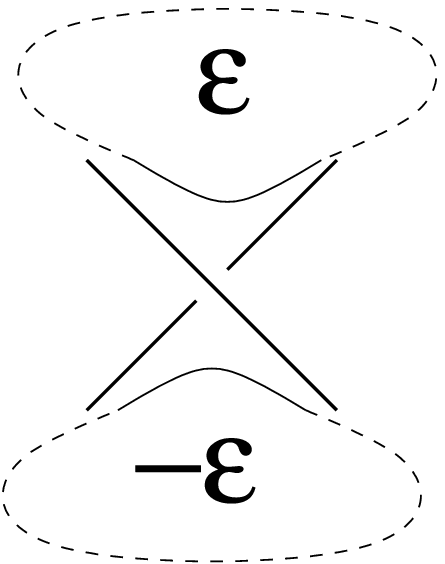,height=1.5cm}}& &
{\psfig{figure=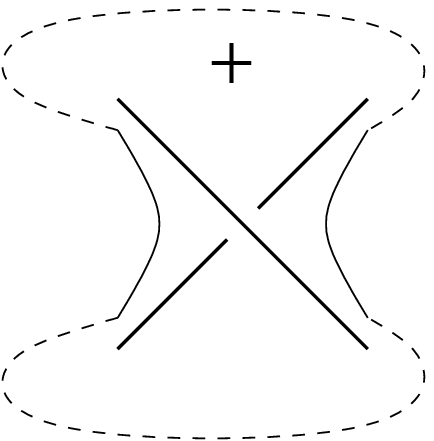,height=1.5cm}}\\ \hline
{\psfig{figure=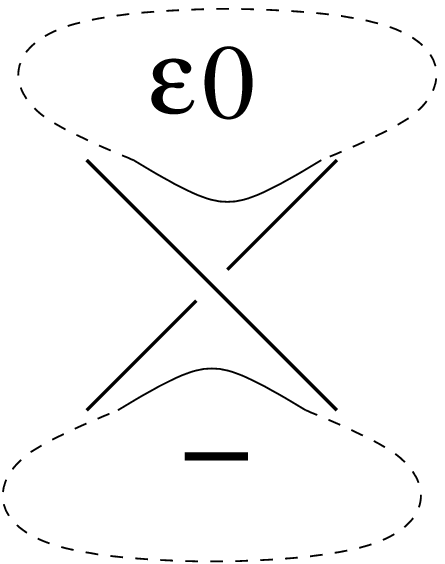,height=1.5cm}}& &
{\psfig{figure=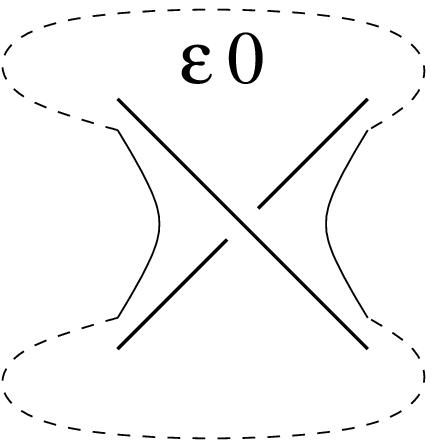,height=1.5cm}}\\ \hline
{\psfig{figure=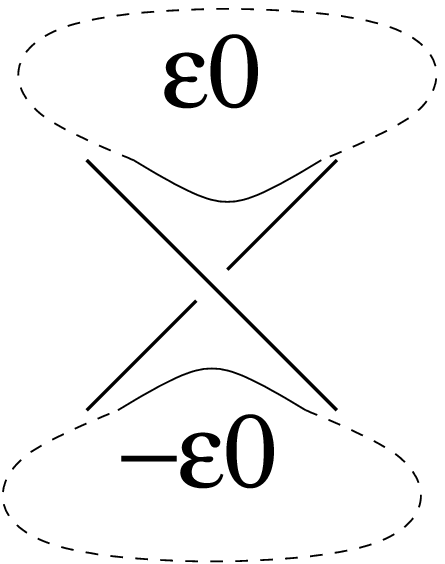,height=1.5cm}} &  &
{\psfig{figure=monep.eps,height=1.5cm}}\\ \hline
{\psfig{figure=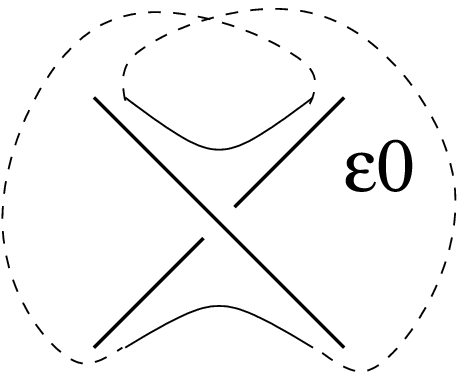,height=1.5cm}}& & 
{\psfig{figure=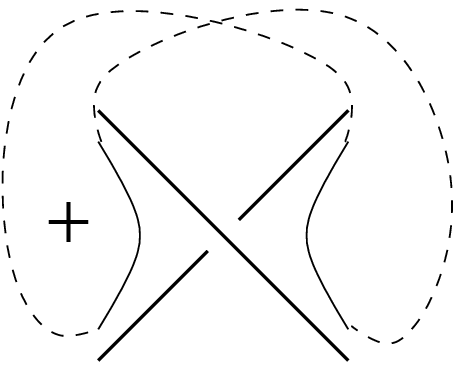,height=1.5cm}}\\ 
\hline
\end{array}$  \end{center} 
\bigskip
\centerline{Table 2.1: Neighborhood of $v$}

Our convention is to label the trivial circles by $+$ or $-$ and the
non-trivial circles by $+0$ or $-0,$ to emphasize that the signs of
the non-trivial circles do not count towards the $j$-grading.
Notice that the states in last row labeled by $\varepsilon 0$ bound
a M\"obius band and, therefore, they may exist in unorientable surfaces only.
The table has the following symmetry: it does not change
by reversing all signs and markers together with flipping $S$ with $S'.$ 

Note also that the degree of $d$ is $-2$ and that $C_{i,*,*}(D)=0$
either for all odd or all even $i,$ depending on the number of 
crossings of $D.$

Finally, notice also that for any link diagram $D$ in $D^2,$ 
$C_{ijs}(D)=0$
for $s\ne 0$ and that $C_{i,j,0}(D)$ and $d$
coincide with Viro's $C_{i,j}(D)$ and $d,$ \cite{V}, 
for framed links.

%
\section{Khovanov homology}
\label{s_d^2}
%

\begin{theorem}[Proof in Section \ref{sec_dd}]\label{d^2}$\phantom{9}$

{\rm(1)}\qua If $F\ne \RP^2$ and $D$ is an unoriented link diagram in $F$ with
ordered crossings, then $d^2=0.$ Hence $(C_{*js}(D), d_{*js})$ is a 
chain complex for any $j\in \Z,$ $s\in \Z\C(F).$

{\rm(2)}\qua For any unoriented link diagram in $\RP^2$ with
ordered crossings, $d^2=0$ mod $2.$ Hence $(C_{*js}(D)\otimes
\Z/2, d_{*js})$ is a chain complex for any $j\in \Z,$ $s\in \Z\C(F).$
\end{theorem}

Theorem \ref{d^2}(1) does not hold for ${\mathbb R}{\mathbb P}^2$.
\begin{example}\label{example}
For the knot diagram $D$ in $\RP^2$ shown below, $C_{2,0,0}(D)$
is generated by the single enhanced state $S$ given by ``$+$'' markers at 
$v$ and $w$ and the single circle of $S$ labeled by ``$-$'',
$$\RP^2= \parbox{1.8cm}{\psfig{figure=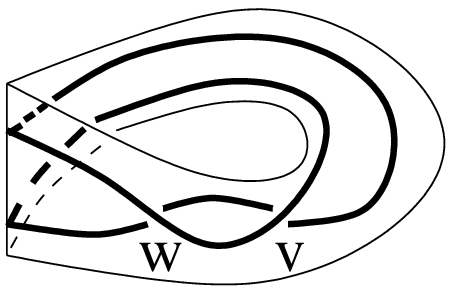,height=1.5cm}}\ 
\hspace*{.2in} \cup_{S^1} D^2.$$
Since $d_w d_v(S)=0$ and $d_v d_w(S)\ne 0,$ $d^2(S)\ne 0.$
\end{example}

From now on we assume that $F\ne \RP^2.$
For any group $G$ we call the homology groups 
of $(C_{*js}(D)\otimes G, d_{*js}),$
$$H_{ijs}(D;G) = {\rm Ker}\ d_{ijs}/{\rm Im}\ d_{i+2,j,s},$$
the Khovanov homology groups of $D$ with coefficients in $G.$ 

For  $G=\Z$ we simplify this notation to $H_{ijs}(D).$

%
\section{Change of crossing order and Reidemeister moves}
\label{main_statements}
%

\begin{theorem}\label{independence} {\rm (Proof in Section 
\ref{sec_independence})}\qua
Given any two orderings of crossings of $D$ and the corresponding
differentials $d_1,d_2$ on $C(D)$ there exists a natural chain
isomorphism $f_{12}: (C(D),d_1)\to (C(D),d_2).$ 
The term ``natural'' means that for any three orderings of crossings
of $D,$ $f_{12}f_{23}=f_{13}.$
\end{theorem} 

\begin{theorem}\label{main_thm}{\rm (Proof in Sections \ref{r_I},
    \ref{r_II} and \ref{r_III})}\qua
Let $D$ be a link diagram in $F.$

{\rm(1)}\qua For any $i,j\in \Z,$ $s\in \Z\C(F),$ 
$H_{ijs}(D)$ is preserved (up to an 
isomorphism) by the Reidemeister moves II and III.
 
{\rm(2)}\qua For a diagram 
\parbox{.8cm}{\psfig{figure=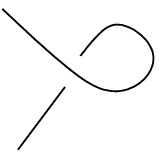,height=0.6cm}}
obtained from \parbox{.8cm}{\psfig{figure=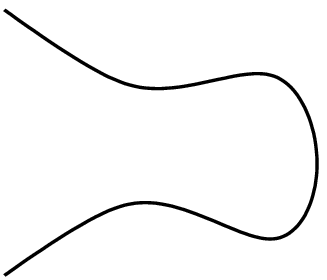,height=0.6cm}}
by adding a negative twist, the map
$\parbox{.8cm}{\psfig{figure=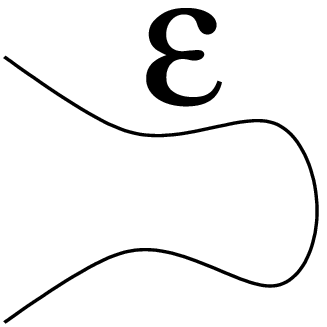,height=.6cm}}\to 
\parbox{.8cm}{\psfig{figure=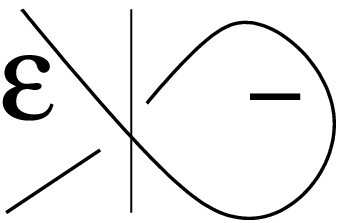,height=.6cm}}$\hspace*{.1in} for 
$\varepsilon=\pm, \pm0,$ induces a chain map $\rho_I:C_{ijs}(
\parbox{.8cm}{\psfig{figure=kinksmooth.eps,height=0.6cm}})\to
C_{i-1, j-3,s}(\parbox{.8cm}{\psfig{figure=kink1.eps,height=0.6cm}})$
which yields an isomorphism of homology groups,
$$\rho_{I*}:H_{ijs}(
\parbox{.8cm}{\psfig{figure=kinksmooth.eps,height=0.6cm}})\to
H_{i-1, j-3,s}(\parbox{.8cm}{\psfig{figure=kink1.eps,height=0.6cm}}).$$ 
\end{theorem}

Consequently, if $M$ is an orientable $I$-bundle over $F$,
then $H_{ijs}(D)$ is an invariant of (band) links in $M$ under 
ambient isotopy.

%
\section{Viro's exact sequence}
\label{s_viro}
%

For any skein triple

\centerline{\psfig{figure=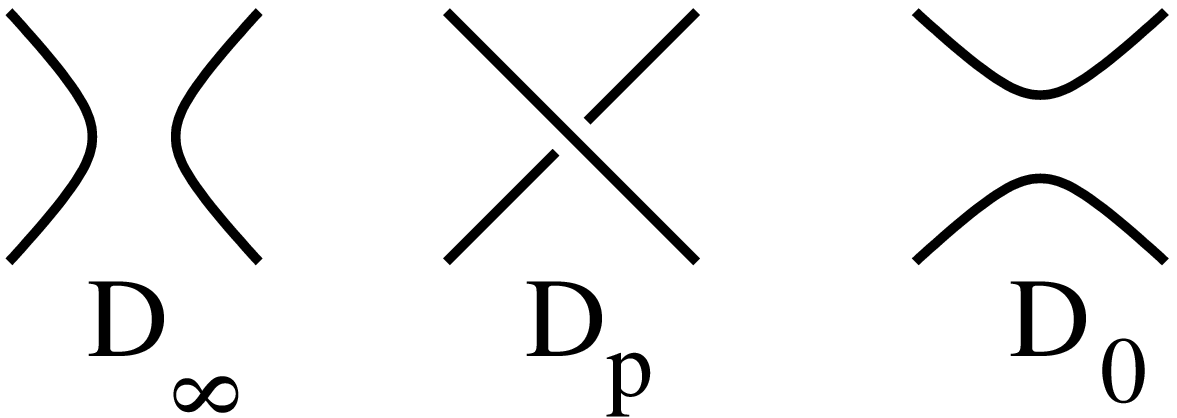,height=1.4cm}}

of link diagrams in a (possibly non-orientable) surface $F$ consider the map 
$\alpha_0: C_{ijs}(D_{\infty}) \to 
C_{i-1,j-1,s}(D_p)$ given by the embedding shown in Figure \ref{s_viro}.1(a) 
and the map $\beta: C_{ijs}(D_p) \to 
C_{i-1,j-1,s}(D_0)$ which is the projection shown in
Figure \ref{s_viro}.1(b).\vspace*{.2in}

\centerline{\psfig{figure=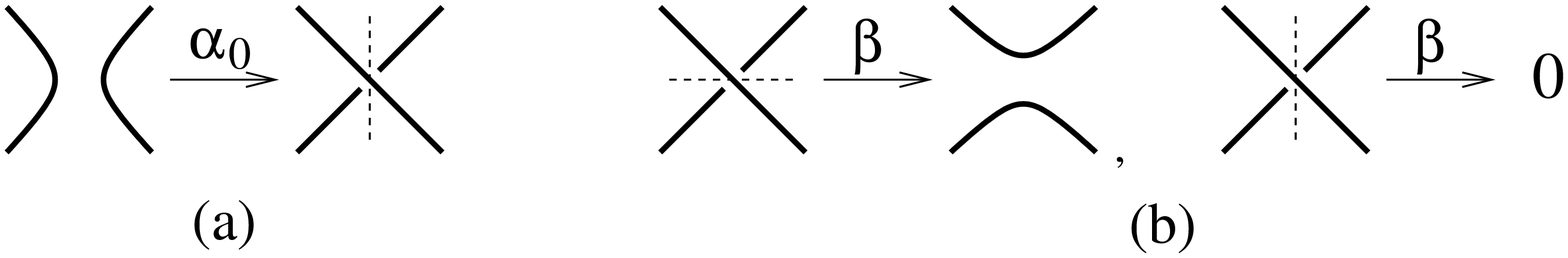,height=1.5cm}}
\begin{center}
{\small Figure \ref{s_viro}.1}
\end{center}

The following is a generalization of the exact sequence 
of \cite[Section 3.3]{V}.

\begin{theorem}\label{ses}
Let $F$ be any surface (including $\RP^2$) and let $D_\infty,$ $D_p,$
$D_0$ be a skein triple of any link diagrams in $F$ as above.
Let the orderings of the crossings in $D_\infty$ and in $D_0$ be
inherited from the ordering of crossings in $D_p$, and let 
$\alpha(S)=(-1)^{t'(S)}\alpha_0(S)$ where
$t'(S)$ is the number of negatively labeled 
crossings in $S$ before $p.$ Then
\begin{enumerate}
\item the maps $\alpha: C_{ijs}(D_{\infty}) \to 
C_{i-1,j-1,s}(D_p),$ $\beta: C_{ijs}(D_p) \to 
C_{i-1,j-1,s}(D_0)$ are chain maps, and
\item the sequence
\begin{equation}\label{eses}
0 \to C(D_{\infty}) \stackrel{\alpha}{\to}
C(D_p) \stackrel{\beta}{\to} C(D_{0}) \to  0
\end{equation}
is exact.
\end{enumerate}
\end{theorem}

\begin{proof}
(1)\qua Note that $\dh_p(\alpha(S))=0.$ (The partial derivative $\dh_p$
was defined in Section \ref{s_diff}.)
Since $[S:S']_q=[\alpha_0(S):\alpha_0(S')]_q$ for any $q\ne p,$
$\alpha$ and $\dh_q$ commute up to sign.
Furthermore,
$$t'(S)+t(\alpha_0(S),q)=t'(S_i)+t(S,q)\ {\rm mod}\ 2,$$
for any summand $S_i$ of $\dh_q(S)$ and, therefore,
$\dh_q(\alpha(S))= \alpha(\dh_q(S)).$
Hence, 
$$d(\alpha(S))=\sum_{q\ne p} {\hat d}_q(\alpha(S))=\alpha(\sum_q {\hat
d}_q(S))= \alpha(d(S)).$$
We have $\beta(d_p(S))=0$ and $d_q(\beta(S))=\beta(d_q(S))$ for
$S\in C(D_p)$ and any crossing $q$ of $D_p.$ 
Since $t(\beta(S),q)=t(S,q)$ for any $q\ne p,$ $\dh_q(\beta(S))=
\beta(\dh_q(S)),$ and, hence, $\beta$ is a chain map as well.

(2) is obvious.
\end{proof}

The short exact sequence (\ref{eses}) leads to the
following long exact sequence of homology groups:
\begin{equation}\label{les}
... \to H_{ijs}(D_{\infty}) \stackrel{\alpha_*}{\to} H_{i-1,j-1,s}(D_p)  
 \stackrel{\beta_*}{\to}
H_{i-2,j-2,s}(D_0) \stackrel{\partial}{\to} H_{i-2,j,s}(D_{\infty})
 \to ...
\end{equation}
Let $\bar\alpha_0: C_{ijs}(D_p) \to C_{i+1,j+1,s}(D_\infty)$
and $\bar\beta: C_{ijs}(D_0) \to C_{i+1,j+1,s}(D_p)$ 
be the following homomorphisms of abelian groups: 
$$\bar\beta(\pic{L0})=(\pic{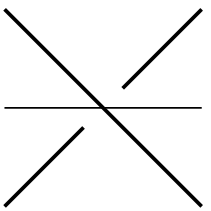})\quad {\rm and\quad } 
\bar\alpha_0(\pic{L+pos})=0,\ \bar\alpha_0(\pic{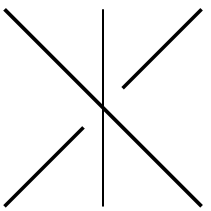})=(\pic{Linfty}).$$
Finally, let $\bar\alpha(S)=(-1)^{t'(S)}\bar\alpha_0(S)$ and let
$\gamma,\hat \gamma: C_{ijs}(D_0)\to C_{i,j+2,s}(D_\infty)$ be
defined as
\begin{equation}\label{gamma}
\gamma=\bar\alpha_0 d_p \bar\beta
\quad {\rm and}\quad \hat \gamma(S)=
\bar\alpha\hat d_p \bar\beta=(-1)^{m(S)}\gamma(S),
\end{equation}
where $m(S)$ denotes the number of negative crossing markers in $S.$

The following properties of $\gamma$ and $\hat \gamma$ will be 
needed later:

\begin{proposition}\label{gamma_prop} Let $F\ne \RP^2.$
\begin{enumerate}
\item $\hat \gamma$ restricted to $Ker\, d\subset 
C_{ijs}(D_0)$ defines the connecting homomorphism $\partial:
H_{ijs}(D_0)\to H_{i,j+2,s}(D_\infty).$

\item $d\hat \gamma=-\hat \gamma d.$
 
\item $d\gamma=\gamma d.$ Consequently, $\gamma:
C_{ijs}(D_0)\to C_{i,j+2,s}(D_\infty)$ 
is a chain map. (Note, however, that neither
$\bar\alpha_0$ nor $\bar\beta$ is a chain map.)

\item For any enhanced state $S\in {\S}_{ijs}(D_0),$
$\gamma(S)=\sum S',$ where the sum is over all enhanced states $S'\in
{\S}_{i,j+2,s}(D_\infty)$ obtained from $S$ by changing the 
\pic{L0}-smoothing
at $p$ to \pic{Linfty}-smoothing in such way that the circles which do 
not pass through the smoothings of $p$ in $S$ and $S'$ have identical labels.
\end{enumerate}
\end{proposition}

\begin{proof}(1)\qua Since $\bar\alpha \hat d_q \bar\beta=0$ for $q\ne p,$
$\hat\gamma=\bar\alpha \hat d_p \bar\beta=\bar\alpha d \bar\beta.$

Since $\bar\alpha$ and $\bar\beta$ split the exact sequence
of groups (\ref{eses}), $\hat \gamma$ restricted to 
$Ker\, d$ defines the connecting differential $\partial.$ 

(2)\qua  Notice that 
\begin{equation}\label{alphabeta_com}
\dh_q(\bar\beta(S))=\bar\beta(\dh_q(S))\quad {\rm and}\quad
\dh_q(\bar\alpha(S))=\bar\alpha(\dh_q(S))\quad {\rm for}\ q\ne p.
\end{equation}
The second equality follows from 
$$t'(S_q)+t(S,q)=t'(S)+t(\bar\alpha_0(S),q)\quad {\rm mod\ } 2,$$ where
$S_q$ is any state such that $[S:S_q]_q=1.$ By equation (\ref{dvdw})
in Section \ref{sec_dd}, $d_p$ and $d_q$ commute and, therefore, 
$\dh_p\dh_q=-\dh_q\dh_p.$ 
Since $\bar\alpha \dh_q \bar\beta=0$ for $q\ne p,$
$$d \hat \gamma(S)=
d(\bar\alpha \dh_p \bar\beta(S))=\sum_{q} 
\dh_q\bar\alpha \dh_p \bar\beta(S) = - \sum_{q} \bar\alpha \dh_p
\bar\beta(\dh_q S)= -\hat \gamma d(S).$$
(3) follows from (2) by assuming that $p$ is the highest (i.e. the last) of the
vertices of $D_p.$

(4) is left to the reader.
\end{proof}

Note that the definitions of $\gamma, \hat \gamma$ do not depend on the
position of $p$ in the ordering of vertices of $D_p.$ In fact,
by Proposition \ref{gamma_prop}(4), these
maps can be defined without referring to $D_p$ at all.

%
\section{Categorification of the skein modules of $I$-bundles over
 surfaces}\label{s_categorify}
%

Khovanov homology groups, $H_{ij}(L),$ for $L\subset \R^3,$ 
normalized according to \cite{V} (for unoriented,
framed links), satisfy the following equation
\begin{equation}\label{maineq}
\chi_A(H_{**}(L))=[L],
\end{equation}
where $[\cdot]$ is the Kauffman bracket normalized by $[\emptyset]=1,$ and
$\chi_A(H_{**})$ is a version of the polynomial Euler characteristic 
defined for the purpose of this paper as
\begin{equation}\label{eech}
\chi_A(H_{**})=\sum_{i,j} A^j (-1)^\frac{j-i}{2} rk H_{ij}.
\end{equation}
Before formulating a version of identity (\ref{maineq}) for links 
in an $I$-bundle $M$ over 
a surface $F$, one must realize that the natural generalization of 
the Kauffman bracket for band
links $L\subset M$ is their representation in the skein module of 
$M,$ $[L]\in \S_b(M;\Z[A^{\pm 1}]).$
The bracket $[L]$ can be identified with a polynomial only if 
$\S_b(M;\lp)=\lp,$ for
example for $L\subset D^2\times I$ and $L\subset S^2\times I.$
Nonetheless, by Theorem \ref{skeinmod}, 
$[L]$ presents uniquely as 
$$[L]=\sum_{b\in {\mathcal B}(F)} p_b(L) b$$
and, therefore, the ``Kauffman bracket'' of $L$ can be thought as the set of 
polynomials $p_b(L)\in \Z[A^{\pm 1}]$ determined by the above
equation. 
We are going to show that polynomial Euler characteristics of 
our homology groups $H_{ijs}(L)$ determine
$p_b(L)$ for all $b\in \B(F),$ if $F$ is orientable. In other words, 
our homology groups ``categorify'' the coefficients
of links expressed in the natural bases of the skein modules of
$F\times I$ for orientable surfaces $F.$
A slightly weaker statement holds also for twisted $I$-bundles
over unorientable surfaces.

Note that the group $G=\prod_{\C(F)} \{\pm 1\}$ acts on $\Z\C(F),$ by
$$(\ve_\gamma) \cdot \sum_{\gamma\in \C(F)} c_\gamma\gamma \to 
\sum_{\gamma\in \C(F)} \ve_\gamma c_\gamma\gamma.$$
Furthermore, it is not difficult to prove that 
\begin{equation}\label{symm}
H_{i,j,gs}(D)=H_{i,j,s}(D), \text{for any $g\in G,s\in \Z\C(F),$ and 
$i,j\in \Z.$}
\end{equation}
Hence it is enough to consider $s=\sum
c_\gamma \gamma,$ with $c_\gamma\geq 0.$
Such elements form the semigroup $\N\C(F),$ where
$\N=\{0,1,...\}.$ Since $\B_{nb}(F)$ embeds into $\N\C(F)$ in a
natural way, one might expect that $p_b(L),$ for $b\in \B_{nb}(F),$
is the polynomial Euler characteristic of $H_{**b}(L),$ but that is
not the case. In order to explicate the relationship between these
polynomials we need to consider the ring of Laurent polynomials
whose set of formal variables $x_\gamma$ is in $1-1$ correspondence
with elements of $\C(F).$ Consider a map 
$$\phi:\{\text{non-trivial simple closed curves in $F$}\} \to
\Z[x_\gamma^{\pm 1}: \gamma\in \C(F)],$$
sending each $\gamma\in \C(F)$ to $x_{\gamma}+x^{-1}_{\gamma},$
and sending every curve $\gamma$ bounding a M\"obius band to $\phi(\gamma)=2.$
This map extends multiplicatively to
$$\phi: \B(F)\to \Z[x_\gamma^{\pm 1}: \gamma\in \C(F)],$$
if $b=\bigcup_{i=1}^d \gamma_i,$ where $\gamma_i$ are (possibly
parallel) non-trivial simple closed curves in $F,$ then
$\phi(b)=\prod_{i=1}^d \phi(\gamma_i).$
Finally, by Theorem \ref{skeinmod}, $\phi$ extends to
$$\phi: \S(M;\Z[A^{\pm 1}])=\Z[A^{\pm 1}]\B(F)\to \Z[A^{\pm 1}, 
x_\gamma^{\pm 1}: \gamma\in \C(F)],$$
in an obvious way.

Hence we have $\phi([L])=\sum q_m(L) m,$ where the sum is
over all monomials\footnote{By a monomial we mean here a product of variables
with leading coefficient $1$.} $m$ in variables $x_\gamma^{\pm 1},$
for $\gamma\in
\C(F),$ and $q_m(L)\in \Z[A^{\pm 1}].$
Note that the set of such monomials can be identified with $\Z\C(F).$

\begin{theorem}\label{polyEuler}
For any $L$ in an orientable $I$-bundle over $F$ and for any $s\in \Z\C(F),$
$\chi_A(H_{**s}(L))=q_s(L).$
\end{theorem}

\proof
By the Kauffman bracket skein relations,
$$\phi\left(\left[\pic{Lp}\right]\right)=A
\phi\left(\left[\pic{L0}\right]\right)+A^{-1}\phi
\left(\left[\pic{Linfty}\right]\right),$$
$$\phi([L\cup \bigcirc])=-(A^2+A^{-2})\phi([L]),$$ 
and, consequently,
$$q_s\left(\pic{Lp}\right)=Aq_s\left(\pic{L0}\right)+A^{-1}q_s
\left(\pic{Linfty}\right),\quad
q_s(L\cup \bigcirc)=-(A^2+A^{-2})q_s(L).$$
Additionally, if $L$ is a link diagram in $F$ and $M\subset F$ is curve
disjoint from $L$ bounding a M\"obius band then $\phi([L\cup M])=2\phi([L]),$
and consequently, $q_s(L\cup M)=2q_s(L).$

On the other hand, (\ref{eech}) and (\ref{les}) imply
$$\chi_A\left(H_{**s}\left(\pic{Lp}\right)\right)=A
\chi_A\left(H_{**s}\left(\pic{L0}\right)\right)+A^{-1}
\chi_A\left(H_{**s}\left(\pic{Linfty}\right)\right).$$
Additionally, since 
$$H_{ijs}(L\cup \bigcirc)= H_{i,j+2,s}(L)\oplus H_{i,j-2,s}(L),$$
we have $$\chi_A(H_{**s}(L\cup \bigcirc))= -(A^2+A^{-2})
\chi_A(H_{**s}(L)).$$
Finally, if $L$ is a link diagram in $F$ and $M\subset F$ is curve
disjoint from $L$ bounding a M\"obius band then
$$H_{ijs}(L\cup M)= H_{ijs}(L)\oplus H_{ijs}(L)$$
and $$\chi_A(H_{**s}(L\cup M))= 2\chi_A(H_{**s}(L)).$$
By above equations, it is enough to prove the
statement for links in $\B_{nb}(F).$ 
Assume hence that $L$ has a diagram composed of $n_i>0$ curves
$\gamma_i$ (not parallel to each other) for
$i=1,...,d.$ Since the diagram of $L$ has no crossings, we have
$H_{ijs}(L)=C_{ijs}(L)$ and $H_{ijs}(L)=0$ for $(i,j)\ne (0,0).$
Consequently, 
\begin{equation}\label{chi}
\chi_A(H_{*,*,s}(L))=rank\, C_{0,0,s}(L).
\end{equation}
If $s=\sum_{i=1}^d k_i \gamma_i$ then
(\ref{chi}) is the number of labelings of components of $L$ by signs $\pm$ such
that the sum of signs of $n_i$ parallel components $\gamma_i$ is
$k_i.$ This is precisely the number of monomials
$x_{\gamma_1}^{k_1}...x_{\gamma_d}^{k_d}$ appearing in the total expansion of
$$\prod_1^d (x_{\gamma_i}+x_{\gamma_i}^{-1})^{n_i}.\eqno{\qed}$$

Since the groups $H_{ijs}(L)$ determine the 
polynomials $q_s(L)$ and $\phi$ is an embedding for an orientable $F,$
we get

\begin{corollary}\label{categorification}
If $F$ is orientable then the homology groups $H_{ijs}(L),$ taken 
over all $s\in \Z\C(F), i,j\in \Z,$ determine the polynomials 
$p_b(L)$ for $b\in \B(F).$
\end{corollary}

More specifically, for any $s=\sum_{i=1}^d k_i\gamma_i$ 
let $\lambda(s)=\prod_{i=1}^d \left(\frac{y_i+\sqrt{y_i^2-4}}{2}\right)^{k_i}.$
\nl
Then $p_s(L)$ is the coefficient of the monomial 
$y_1^{k_1}...y_d^{k_d}$ in the expansion of
$\sum_s q_s(L)\lambda(s).$ 

If $F$ is unorientable, then $\phi$ is not an embedding and the
statement of Corollary \ref{categorification} does not hold.
In this case, for any $b\in \B_{nb}(F),$ let $\B(F;b)$ denote the set of
all basis elements of $F$ obtained by adding to $b$ disjoint closed
curves bounding M\"obius bands. 
We have
$$\B(F)=\amalg_{b\in \B_{nb}(F)} \B(F;b)$$
and the following generalization of Corollary \ref{categorification}
holds:

\begin{corollary}
For any orientable $I$-bundle over $F$ the homology groups 
$H_{ijs}(L),$ taken over all $s\in \Z\C(F), i,j\in \Z,$ 
determine the sums
$$\sum_{b'\in \B(F;b)} 2^{|b'|-|b|} p_{b'}(L)$$ 
for any $b\in \B_{nb}(F).$ Here, $|b|$ denotes the number of components of $b.$
\end{corollary}

%
\section{Cohomology}
\label{cohomology}
%

For a link diagram $D$ in $F,$ the {\em Khovanov cohomology groups}
$H_{ijs}(D)$ are the cohomology groups of the cochain complex 
$(C^{***}(D),d^*),$ for $C^{***}(D)={\rm Hom}(C_{***}(D), \Z).$  
Here is another way of defining the cochain 
complex of $D:$

\begin{proposition}\label{first_isom}
The following diagram commutes:
$$\begin{array}{cc}
C_{ijs}(D)\stackrel{{\tilde d}}{\longrightarrow}& C_{i+2,j,s}(D)\\
\downarrow \psi & \downarrow \psi \\
C^{ijs}(D) \stackrel{d^*}{\longrightarrow} & C^{i+2,j,s}(D),
\end{array}$$
where ${\tilde d}$ is defined by 
$${\tilde d}(S)=\sum_{{\rm all\ states\ } S'\atop  
{\rm all\ crossings\ } v {\rm\ in}\ D} \hspace*{-.2in} 
(-1)^{t(S,v)} [S':S]_v S',$$
and $\psi$ is the isomorphism given by the identification of the basis
composed of enhanced states with its dual. 
\end{proposition}

\begin{proof}
Denote the dual of $x\in C(D)$ by $x^*\in C^*(D).$
We need to prove that
$({\tilde d}(S))^*=d^*(S^*),$ for any $S\in C(D).$
Equivalently, that for any $S_0\in C(D),$
$$({\tilde d}(S))^*(S_0)=d^*(S^*)(S_0).$$
The left side of this equation expands to
$$\sum_{{\rm all\ states\ } S'\atop  
{\rm all\ crossings\ } v {\rm\ in}\ D} \hspace*{-.2in} 
(-1)^{t(S,v)} [S':S]_v (S')^*(S_0) =(-1)^{t(S,v_0)}[S_0:S]_{v_0},$$
where $v_0$ is the crossing of $D$ at which the markers
of $S_0$ and $S$ differ. If $[S_0:S]_v=1$ for some $v$ then such 
crossing $v$ is unique
and $v_0=v.$ On the other hand, we have
$$d^*(S^*)(S_0)=S^*(d(S_0))=(-1)^{t(S_0,v_0)}[S_0:S]_{v_0}.$$
Since the markers of $S$ and $S_0$ differ at $v_0$ only,
$t(S_0,v_0)=t(S,v_0)$ and the proof is complete.
\end{proof}

For a link diagram $D$ we denote its mirror image, i.e. the diagram 
obtained by switching all crossings of $D$, by ${\bar D}.$
In order to relate the cohomology of $D$ to homology of $\bar D,$ 
consider $d^+(S)= \sum_{v} (-1)^{t^+(S,v)} d_{ijs;v}(S),$
where $t^+(S,v)$ is the number of positive markers
assigned to crossings in $S$ bigger than $v.$ 
We can show that $d^+d^+(S)=0$ in the same way as $d^2(S)=0$ --
both statements follow from equation (\ref{dvdw}) in Section \ref{sec_dd}.
Therefore either of the differentials $d,d^+$ could be used for the
definition of the Khovanov chain complex of links.
The following lemma shows that these chain complexes are isomorphic.

\begin{proposition}\label{d+} Let $D$ be a link diagram with $n$ crossings and
let for any enhanced state $S$ of $D,$ $g(S)=(-1)^{u(S)}S,$ 
where $u(S)$ is the number of positively marked crossings $v_i$ in 
$S$ such that $i=n+1$ mod $2.$
This map induces an isomorphism of chain complexes
$g:(C(D),d)\to (C(D),d^+).$
\end{proposition}

\begin{proof}
Both $d(S)$ and $d^+(S)$ are linear combinations of states $S'$
obtained from $S$ by changing one positive marker to a negative one.
Therefore it is enough to prove that if $S,S'$ differ by a marker
at an $i$th vertex $v_i$ then the following morphisms form a commutative
diagram:
$$\begin{array}{ccc}
 & S\to (-1)^{t(S,v_i)}S& \\
S & & S \\
\downarrow &  & \downarrow \\
(-1)^{u(S)}S & & (-1)^{u(S)} S\\
 & S\to (-1)^{t^+(S,v_i)}S & 
\end{array}$$
Since $t^+(S,v_i)=n-i-t(S,v_i)$ and $u(S)-u(S')=n-i$ mod $2$
the above diagram commutes.
\end{proof}

For a crossing $v$ in $D,$ we denote by $\bar v$ the corresponding crossing
in $\bar D.$ We order the crossings of $D$ and of $\bar D$ in a consistent
manner, i.e. if $v$ is the $i$th crossing in $D$ then 
$\bar v$ is the $i$th crossing in $\bar D.$
Now, consider the map
$$\varphi: C_{ijs}(D)\to C_{-i,-j,-s}({\bar D})$$ 
such that for any enhanced state $S$ of $D,$ $\varphi(S)$ is an enhanced
state of $\bar D$ whose every crossing $\bar v$ is labeled by the
marker opposite to the marker of $v$ in $S.$
Under these assumptions the loops of $S$ and $\varphi(S)$
are naturally identified. We assume that $\varphi(S)$ reverses the signs of 
all loops in $S.$

\begin{proposition}\label{second_isom}
The following diagram commutes
$$\begin{array}{ccc}
C_{ijs}(D) & \stackrel{{\tilde d}}{\longrightarrow}&
    C_{i+2,j,s}(D)\\
\downarrow \varphi & & \downarrow \varphi \\
C_{-i,-j,-s}({\bar D})& \stackrel{d^+}{\longrightarrow} &
C_{-i-2,-j,-s}({\bar D}).
\end{array}$$
\end{proposition}

\begin{proof} We need to prove that 
\begin{equation}\label{eq1}
\varphi({\tilde d}(S))=d^+(\varphi(S)), {\rm \ for\ any\ } S\in C(D).
\end{equation}
The left hand side equals
$$\sum_{{\rm all\ } S'\in {\S}(D)\atop  
{\rm all\ crossings\ } v {\rm\ in}\ D} \hspace*{-.2in} 
(-1)^{t(S,v)} [S':S]_v \varphi(S')$$
and the right hand of (\ref{eq1}) is 
$$\sum_{{\rm all\ } {\bar S}\in {\S}({\bar D})\atop  
{\rm all\ crossings\ } {\bar v} {\rm\ in}\ {\bar D}} \hspace*{-.2in} 
(-1)^{t^+(\varphi(S),{\bar v})} [\varphi(S):{\bar S}]_{\bar v} {\bar S}.$$
Since the states ${\bar S}$ of $\bar D$ are in 1-1 correspondence with
states of $D$ via the map $\varphi,$ this expression takes the form
$$\sum_{{\rm all\ } S'\in {\S}(D)\atop  
{\rm all\ crossings\ } v {\rm\ in}\ D} \hspace*{-.2in} 
(-1)^{t^+(\varphi(S),{\bar v})} [\varphi(S):\varphi(S')]_{\bar v} 
\varphi(S').$$
Now the statement follows from the following two identities:
\begin{equation}
t(S,v)=t^+(\varphi(S),{\bar v})
\end{equation}
\begin{equation}
[S:T]_v =[\varphi(T): \varphi(S)]_{\bar v} {\rm \ for\ any\ } 
S,T \in {\S}(D).
\end{equation}
The first one follows from the fact that $\varphi$ reverses the signs of 
markers, and the second one follows directly from the definition of 
the incidence number.
\end{proof}

By combining the statements of Propositions \ref{first_isom}, \ref{d+}, and 
\ref{second_isom}, we conclude 

\begin{corollary} For any diagram $D$ and $i,j\in \Z,$ $s\in \Z\C(F),$
$$g_*\varphi_*\psi_*^{-1}:H^{ijs}(D) \to H_{-i,-j,-s}({\bar D})$$
is an isomorphism.
\end{corollary}

By Universal Coefficient Theorem (cf.\ \cite[Corollary 3.5]{Ha}), 
$$ H_{-i,-j,-s}({\bar D})= H^{ijs}(D) = H_{ijs}(D)/T_{ijs}(D) \oplus 
T_{i-2,j,s}(D),$$
where $T_{ijs}(D)$ is the torsion part of $H_{ijs}(D).$ 
Consequently, the torsions of homology
groups of $D$ and $\bar D$ coincide (up to a shift of indices).
This observation is used in \cite{AP}.

%
%

\section{Khovanov homology for different $I$-bundle structures}

\label{I-bundle}
%
%

Any handlebody $H$ is an $I$-bundle over a surface, although the
$I$-bundle structure of $H$ is not unique. For example, if
\begin{center}
$F_1=\parbox{2cm}{\psfig{figure=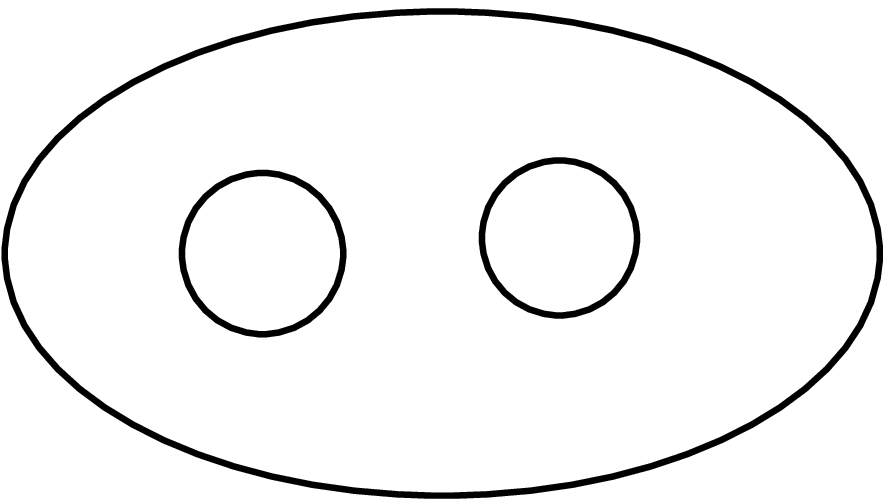,height=1cm}}$ and
$F_2=\parbox{2.1cm}{\psfig{figure=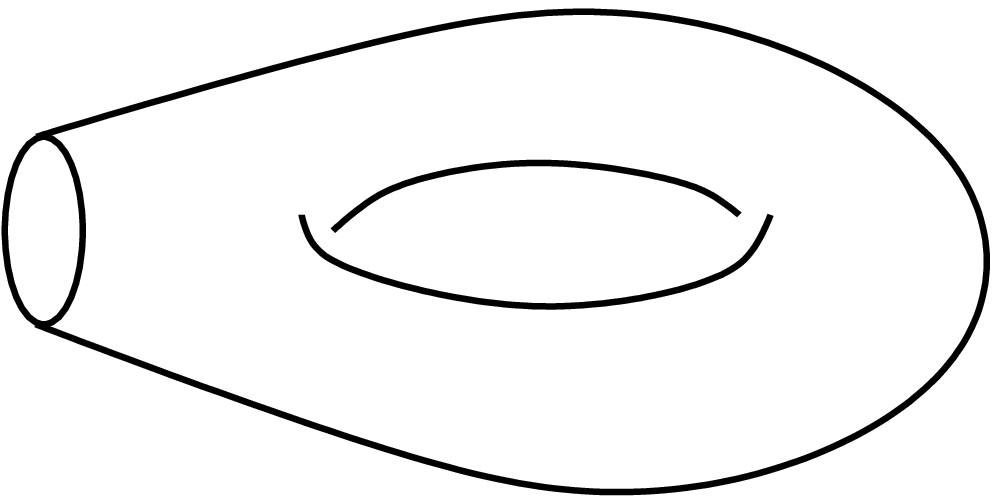,height=1cm}},$
\end{center}
then $F_1\times I\simeq F_2\times I$ are both handlebodies of genus $2.$
Hence, it is natural to
ask whether there exist homology groups of links $L\subset H$ which
do not depend on the $I$-bundle structure of $H$ and
on the choice of $s\in \Z\C(F)$ for the base surface $F.$
One obvious choice to consider is
$${\mathcal H}_{ij}(L)=\bigoplus_{s\in \Z\C(F)} H_{ijs}(L).$$
We are going to give two examples showing that unfortunately 
${\mathcal H}_{ij}(L)$
does depend on the choice of $I$-bundle structure of $H.$

{\bf Example 1}\qua
The spaces $F_2\times I$ and $F_1\times I$ are homeomorphic 
by a homeomorphism mapping the curve $\gamma$ in $F_2\times \{0\}$ to 
the curve $\gamma'$ in $\partial(F_1\times I).$
\begin{center}
{\psfig{figure=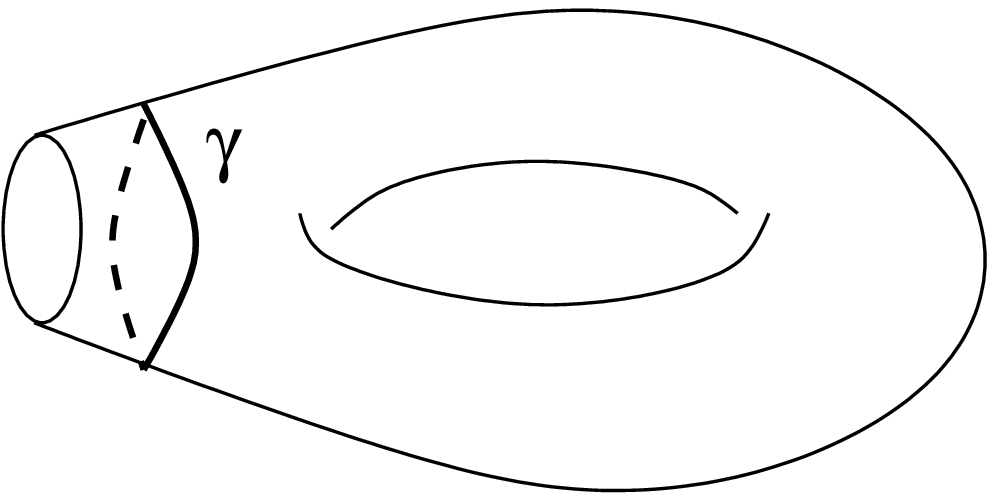,height=1.4cm}} \quad 
{\psfig{figure=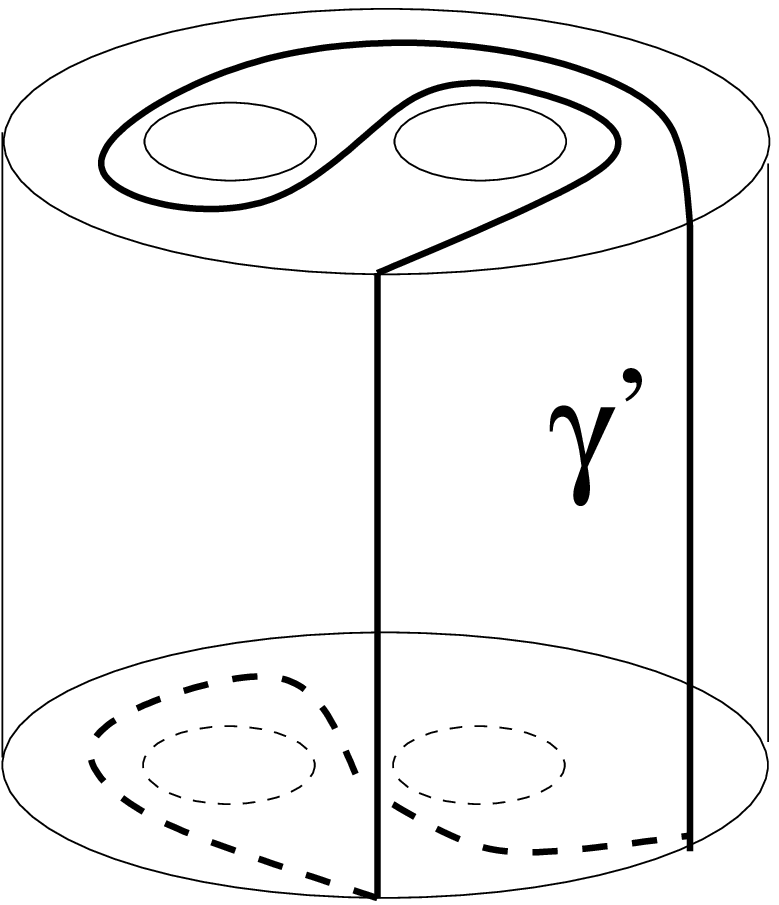,height=2.2cm}}
\end{center}
Therefore, the knot $K$ represented by the diagram $D=\gamma$ in $F_2$
is given by the following diagram $D'$ in $F_1,$

\centerline{{\psfig{figure=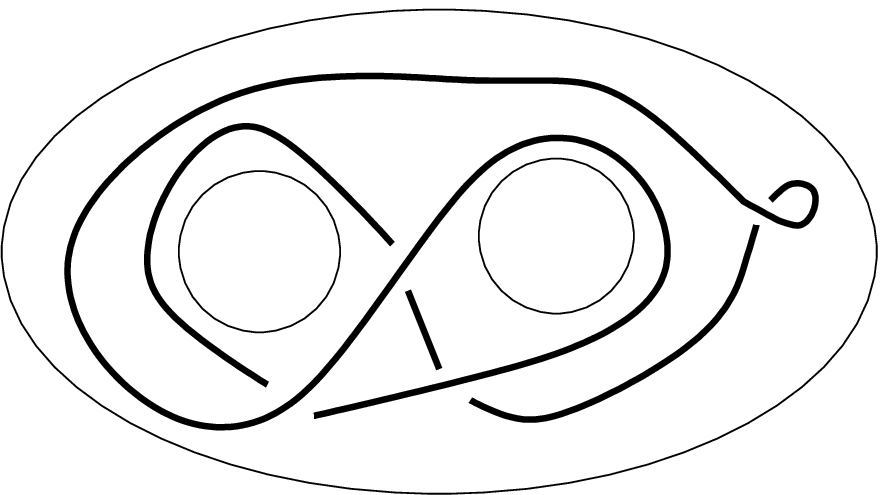,height=2cm}}}

The calculations of homology groups of $D$ and $D'$ yield
$${\mathcal H}_{ij}(D)=\begin{cases} \Z^2 & \text{for $(i,j)=(0,0)$}\\ 0 &
  otherwise\\ \end{cases}$$
$${\mathcal H}_{ij}(D')=\begin{cases} \Z^7 & \text{for $(i,j)=(0,-2)$}\\
\Z^8 & \text{for $(i,j)=(-2,-4)$}\\ \Z^3 & \text{for $(i,j)=(-4,-6)$}\\ 
0 & \text{otherwise.}\\ \end{cases}$$

{\bf Example 2}\qua The homeomorphism of $H_2=F_1\times I$ given by
the following Dehn twist
$${\psfig{figure=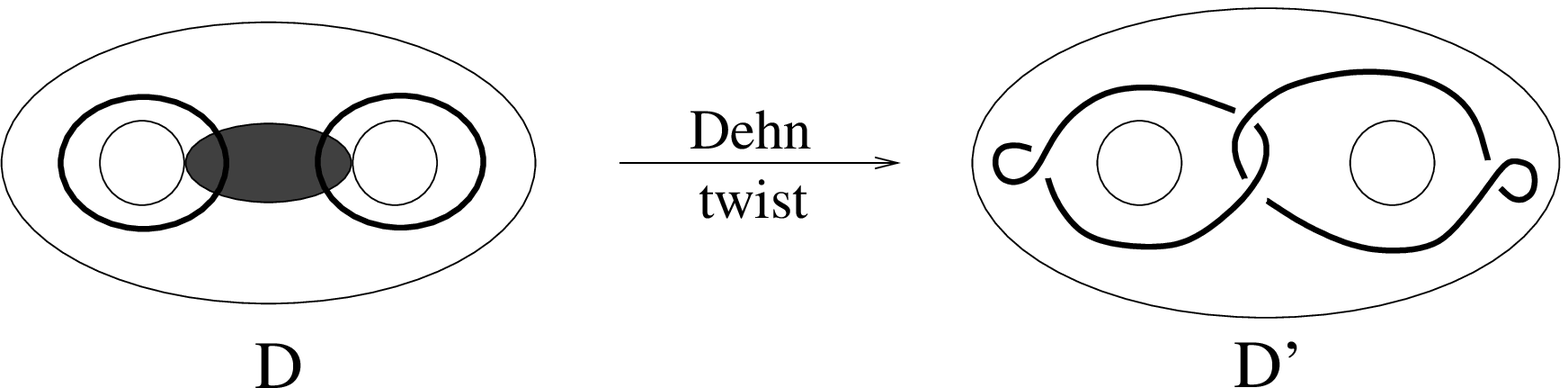,height=1.9cm}}$$
maps $D$ to $D'.$

We have
$${\mathcal H}_{ij}(D)=\begin{cases} \Z^4 & \text{for $(i,j)=(0,0)$}\\ 0 &
  \text{otherwise}\\ \end{cases}\quad\quad {\rm and}$$
$${\mathcal H}_{ij}(D')=\begin{cases} \Z^4 & \text{for $(i,j)=(3,5)$}\\ 
\Z^2 & \text{for $(i,j)=(1,3)$ and $(-1,-1)$}\\ 
0 & \text{otherwise.}\\ \end{cases}$$

%
\section{The proofs}
%
%

%
\subsection{Invariance of homology under crossing order changes}
\label{sec_independence}
%


Consider two different orderings of the vertices of a link diagram $D.$
Let $S$ be an enhanced state of $D.$ Each of the orderings induces
an ordering on the set of negatively marked crossings in $S.$
These two orderings (restricted to negative crossings) differ by a
permutation $\sigma_S.$ We define $f_{12}(S)=sgn(\sigma_S)S$
and extend it to $f_{12}:C(D)\to C(D).$ 

For any three orderings of crossings of $D,$ $f_{12}f_{23}=f_{13}$
and therefore for the remainder of the proof we assume that
the orderings of crossings of $D$ differ by a transposition $(i,i+1).$
We need to prove that $f_{12}(d_1(S))=d_2(f_{12}(S)),$ and for that it
is enough to show that if $[S:S']_v=1$ then
$$(-1)^{t_1(v,S)-t_2(v,S)}=sgn(\sigma)sgn(\sigma'),$$
where $t_1(v,S),t_2(v,S)$ denote the numbers of  
negative markers assigned to vertices above $v$ in $S,$ with respect
to the two orderings. Since $[S:S']_v=1,$ $v$ has the $+$ marker in $S.$
To complete the proof 
note that each of the sides of this equation is equal $-1$ if and only if
either \begin{itemize}
\item $v$ is the $i$th crossing in $D$ (with respect to the first
ordering) and $S$ assigns $-$ to the $(i+1)$st
crossing; or
\item $v$ is the $(i+1)$th crossing in $D$ and $S$
assigns $-$ to the $i$th crossing.
\end{itemize}

%
%
\subsection{Proof of $d^2=0$ (Theorem \ref{d^2})}

\label{sec_dd}
%

For any link diagram $D$ in $F,$ any 
$i,j\in \Z,$ $s\in \Z\C(F),$ and any enhanced 
state $S\in C_{i,j,s}(D),$ we are going to prove that 
$d^2_{*js}(S)=0,$ if $F\ne \RP^2,$ and 
$d^2_{*js}(S)=0$ mod $2$ for any $F.$ The signs of partial derivatives
in $d_{ijs}$ are defined in such
a way that it is sufficient to show the following two equations
\begin{equation}
\label{dvdw}
d_v \circ d_w(S)=d_w \circ d_v(S)\quad {\rm if\ } F\ne \RP^2
\end{equation}
\begin{equation}
\label{dvdw2}
d_v \circ d_w(S)=d_w \circ d_v(S)\quad {\rm mod\ } 2 {\rm\ (for\ any\ } F{\rm)}
\end{equation}
for any crossings $v$ and $w$ of $D$. 

Let $D'$ be a diagram obtained from $D$ by 
smoothing all the crossings of $D$ according to $S$, 
except for $v$ and $w$. Since $S$ can be considered as a state 
of $D'$, it is sufficient to show equation (\ref{dvdw}) for $D'$. 
Therefore, from now on we assume that $D$ has $2$ crossings only: $v$ and
$w.$

\begin{lemma}\label{duality}
If (\ref{dvdw}),(\ref{dvdw2}) hold for all states of $D$ then
they hold for all states of $\bar D$ as well.
\end{lemma}

\begin{proof} If $d^2(S)=0$ for all $S\in C_{2,*,*}(D),$ then
$(d^*)^2(S^*)=0$ for all $S^*\in C^{-2,*,*}(D).$ Hence, by Theorems
\ref{first_isom} and \ref{second_isom},
$d^2(S)=0$ for all $S\in C_{2,*,*}({\bar D}).$
\end{proof}

We denote the trivial circles in any state by T and
the nontrivial circles by N. Hence NT denotes two circles, one trivial
and one nontrivial.
Our proof of (\ref{dvdw}) and (\ref{dvdw2}) does not use the fact that
$D$ lies in a surface and it relies entirely on the following
principle: If $F\ne \RP^2$ then the change of smoothing of any crossing 
of $D$ results in one of the following transformations: 

\begin{center}T $\to$ N, TT or NN, \quad N $\to$ T, N, NT or NN
\end{center}

\begin{center} \label{transf}
TT $\to$ T, \quad NN $\to$ T or N,\quad TN $\to$ N.
\end{center}
Additionally, if $F=\RP^2$ then another possible
transformation is
\begin{equation} \label{TT}
{\rm T} \quad \to \quad {\rm T}.
\end{equation}
(This transformation may not happen in $F\ne \RP^2.$ Proof: A trivial circle
bounds a disk, which after the change of smoothing becomes a M\"obius
band. If the boundary of M\"obius band in $F$ is trivial then $F$ is a union
of a M\"obius band and a disk.)

Since our proof does not use any information about the position of $D$ in $F,$
we can think of any $2$ crossing link diagram in $F$ as an
abstract $4$-valent graph with $2$ rigid vertices.
Any such graph is obtained by connecting the $1$-valent vertices in
the picture below:

\centerline{\psfig{figure=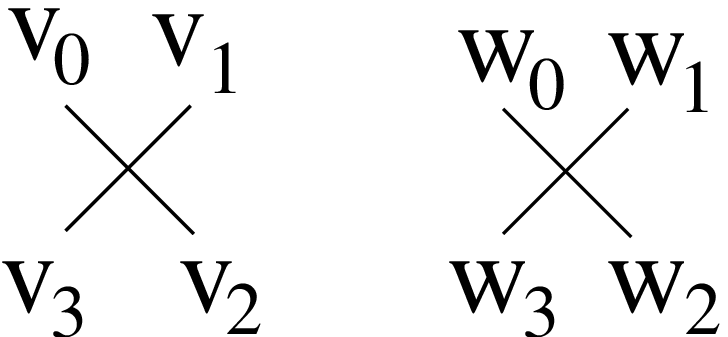,height=.5in}} 

Hence we consider all possible
pairings of the following eight points $v_0,$ $v_1,$ $v_2,$ $v_3,$
$w_0,$ $w_1,$ $w_2,$ $w_3.$
Since obviously (\ref{dvdw}), (\ref{dvdw2}) hold for diagrams $D$ in
which $v$ and $w$ belong to different connected components, without
loss of generality we assume that $v_0,$ $w_0$ are connected by an arc.
The following operations identify isomorphic graphs:
\begin{enumerate} 
\item exchanging $v$ with $w;$
\item  cyclic permutation 
$v_i\to v_{i-k {\rm\ mod\ }4}$ and $w_i\to w_{i-k {\rm\ mod\ }4}$ for
all $i;$ Since $v_0$ and $w_0$ are always assumed to be connected, this 
operation is allowed only if $v_k$ and $w_k$ are connected.
\item flip: $v_i$ is exchanged with $v_{i+2}$ for some $i\in
\{0,1,2,3\}$ while $w_{i+1{\rm\ mod\ }4}$ and $w_{i+3{\rm\ mod\ }4}$
are fixed.  An analogous flip is allowed for vertex $w$ as well.
\end{enumerate}
The last relation follows from the fact that if link diagrams
$D_1,D_2$ are represented by graphs $\Gamma_1,\Gamma_2$ related by a
flip, then there is a natural $1-1$ correspondence between the states
of $D_1$ and the states of $D_2.$ This correspondence is preserved by 
changing of smoothing, incidence numbers, and partial derivatives.

The proof of the following lemma is left to the reader:

\begin{lemma} 
Up to the above relations, there are precisely $5$ different 
abstract connected $4$-valent graphs with $2$ rigid vertices:

\centerline{\psfig{figure=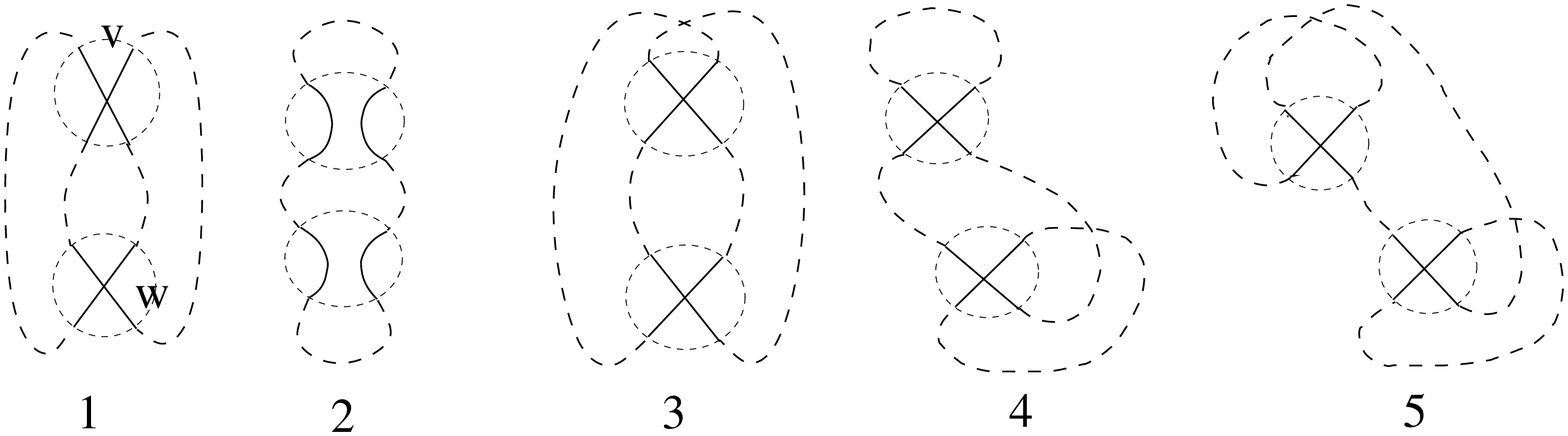,height=1.3in}} 
\end{lemma}

We prove (\ref{dvdw}) and (\ref{dvdw2}) for the above graphs
on a case by case basis. In the proofs we will frequently refer to the 
following arguments for $S:$
\begin{description}
\item[\rm(NZ)] We assume that either $d_vd_w(S)$ or $d_wd_v(S)$ is non-zero.
(Otherwise, clearly (\ref{dvdw}),(\ref{dvdw2}) holds.)
\item[\rm($\tau$)] If $d_v(S)=\sum_i S_i,$ $d_w(S_i)=\sum_j S_{ij}$ then $
\tau(S_i)=\tau(S)+1,\tau(S_{ij})=\tau(S)+2$ for all $i,j.$
(Analogous statement for $v$ and $w$ exchanged).
\item[\rm(S)] Argument by the symmetry of the picture.
\end{description}

\noindent {\bf Notation}\qua
In the proof we use $\underbrace{T...T}_k\underbrace{N...N}_l$ to
denote any state composed of $k$ trivial and $l$ nontrivial
circles. When labels of circles are of importance we use
the symbols \cp,\cm,\ce\ to denote trivial components labeled
by $+,-,$ and $\ve,$ respectively; $\ve$ denotes $+$ or $-$.
The symbols \cpo, \cmo are used to denote the nontrivial circles,
labeled by $+$ and $-,$ respectively.
By abuse of notation, we sometimes write S=TN or S=\cm\cpo
despite the fact that neither TN and \cm\cpo\ denotes a unique state.
In particular, they do not specify which of the circles of S is trivial.
If the order of circles of $S$ is important and if the circles of $S$ are
positioned one above another, then we use
${\cm\atop \cpo}$ to specify that the circle above is trivial and
labeled by $-$ and the bottom one is non-trivial and labeled by $+.$

\centerline{\psfig{figure=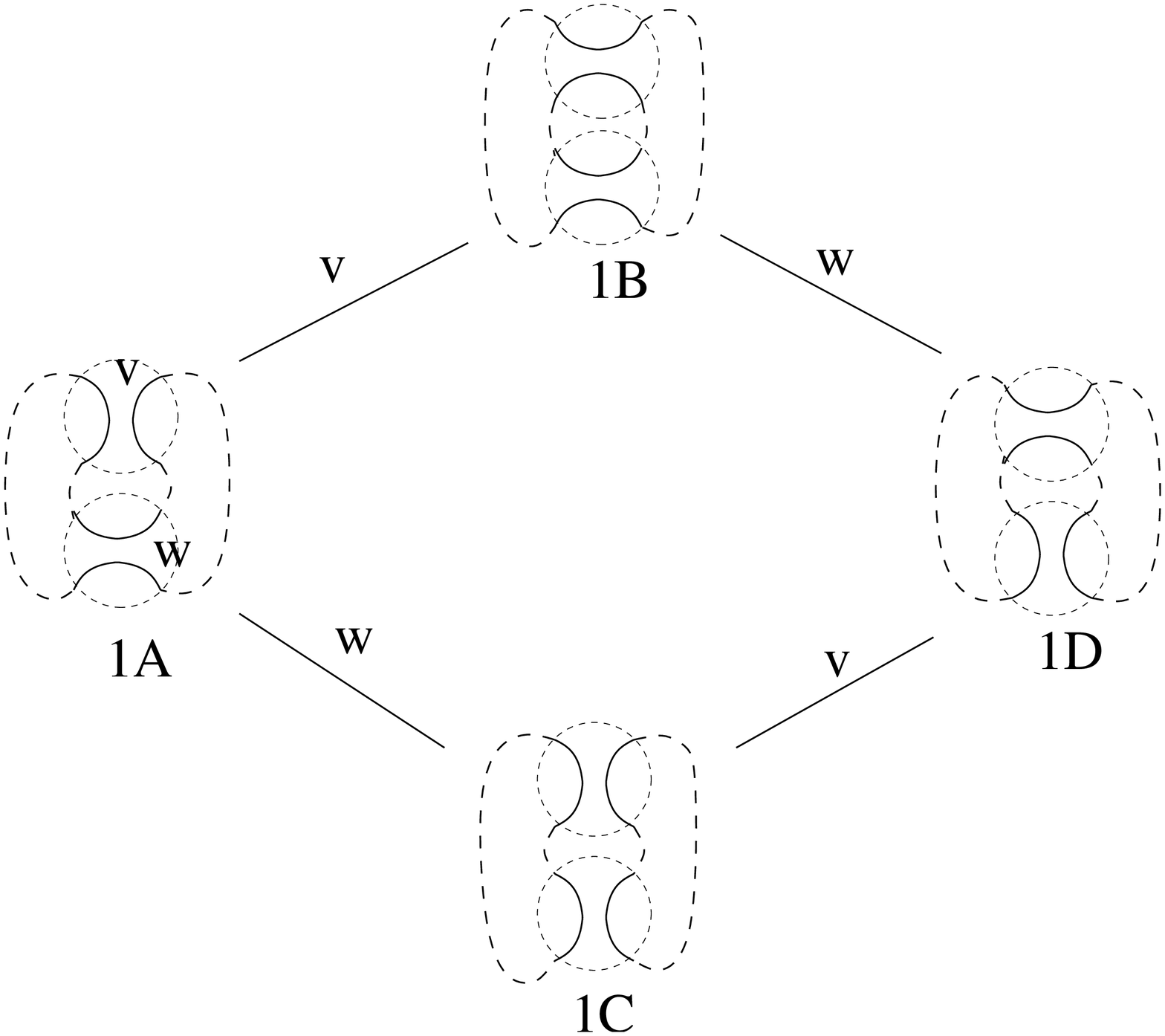,height=2in}}  

We are going to prove (\ref{dvdw},\ref{dvdw2}) for states of the above
four types, case by case.
In each of these cases, the edges with appropriate orientation  will
represent the partial derivatives. For example, if $S$ if of type 1A
then the edges are oriented outwards 1A. Notice that $d_v(S)$
is a sum of states of type 1B, etc.

\begin{description}
\item[\rm(1A)]
By (NZ,\t), 1A=\cm,1D=\cp. Hence 1A is a trivial loop, and
therefore each of the states 1B and 1C is either composed of two 
trivial circles or two nontrivial circles. In either case 
$d_vd_w(\cm)=d_wd_v(\cm)=2\cp.$

\item[\rm(1B)] If 1A and 1B are either both trivial or
both nontrivial then by (S), (\ref{dvdw}) holds. Hence, again by (S),
we assume that 1A is T and 1D is N, and therefore 1B and 1C are NN.
Now by (\t), both sides of (\ref{dvdw}) are $0.$

\item [\rm(1C) and (1D)] The proof follows from (1A) and (1B) by Lemma
\ref{duality}. (Alternatively, the proof follows here from symmetry of
the above diagram.)
\end{description}

\centerline{\psfig{figure=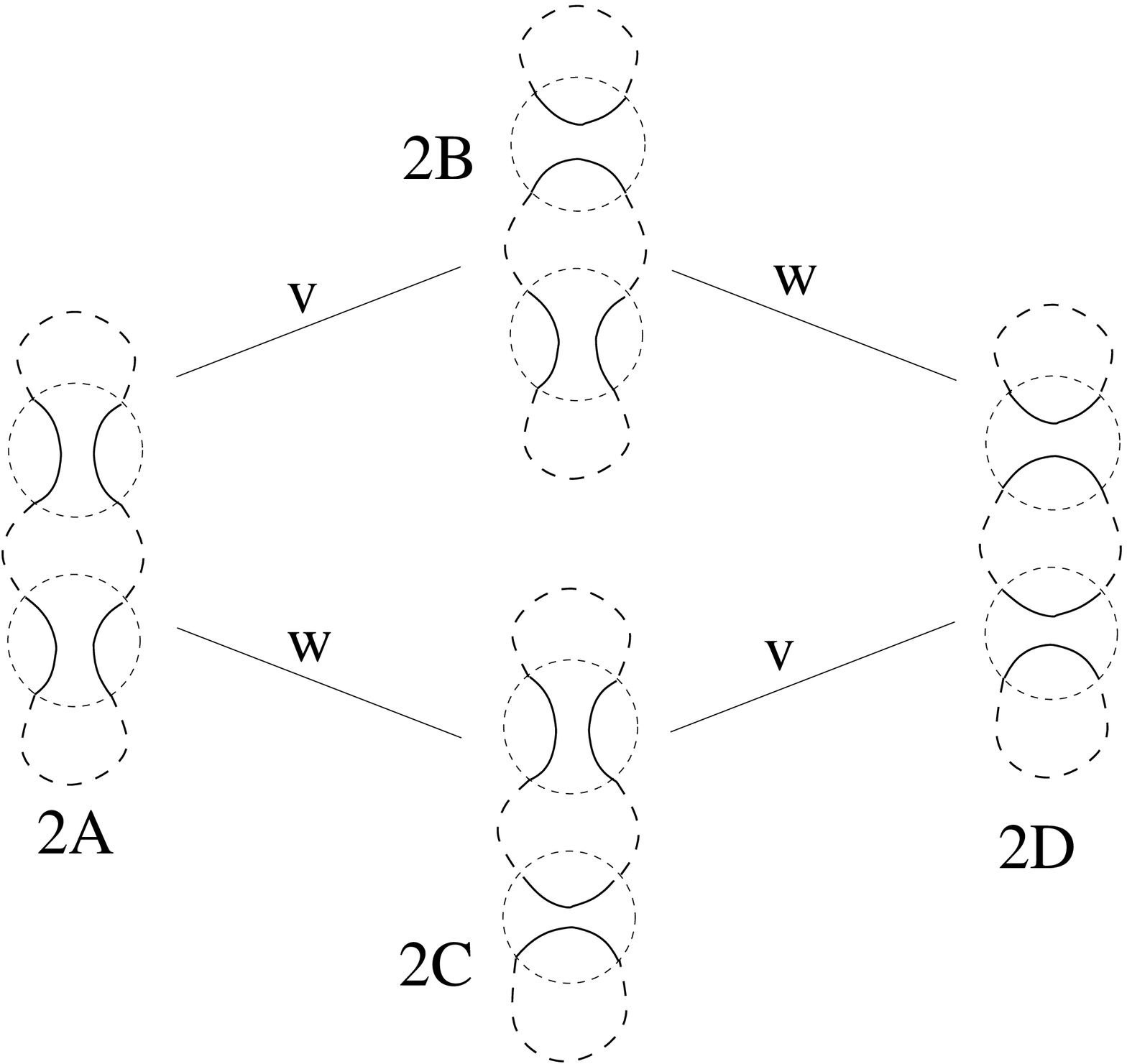,height=2in}}  

\begin{description}
\item[\rm(2A)] If 2D is TTT then 2B,2C are TT and the statement follows
by (S). Hence assume that 2D is not TTT. 
Furthermore, by (\t) and (NZ) we assume that if 2A is N then 2D is TTN.
Therefore, we have the following possibilities:

$\diagtwoa{T}{T}{T}{T}{N}{N}{N}{N}$
\begin{minipage}{4in}
By (\t), $d_vd_w(S)=d_wd_v(S)=0,$ unless $S=\cm.$

On the other hand,
$d_vd_w(\cm)=d_wd_v(\cm)=\three{\cp}{\cpo}{\cmo}+\three{\cp}{\cmo}{\cpo}.$
\end{minipage}

$For \diagtwoa{T}{T}{N}{T}{T}{N}{T}{N},$
$d_vd_w(\ceo)=d_wd_v(\ceo)=\three{\cp}{\cp}{\ceo}.$

\begin{minipage}{2in}
The remaining cases\nl follow by symmetry:
\end{minipage}
$\diagtwoa{T}{N}{N}{N}{T\, {\rm or\, } N}{N}{N}{N},$
$\diagtwoa{N}{T}{N}{T}{N}{T}{N}{T}$ 


\item[\rm(2B)] By (\t, NZ) we exclude the cases
$(2A,2C)=(NN,NN),(NN,TN),$ $(TN,NN).$ The remaining cases are:

$\diagtwob{N}{T}{N}{T}{T}{N}{T}{N},$
$\diagtwob{N}{T}{N}{T}{N}{T}{N}{T},$ 
$\diagtwob{N}{N}{T}{N}{N}{N}{T}{N}.$

By (\t), in each of these cases, $d_vd_w(S)=d_wd_v(S)=0,$ unless
$S=$\cm\ceo.  For $S=$\cm\ceo, we get
$d_vd_w(S)=d_wd_v(S)={\cp\atop \ceo}.$

$\diagtwob{T}{T}{T}{T}{N}{N}{N}{N}$
\begin{minipage}{4in} By (\t), $d_vd_w(S)=d_wd_v(S)=0,$ unless
$S={\cm\atop \cm}.$

For $S={\cm\atop \cm},$ 
$d_vd_w(S)=d_wd_v(S)={\cpo\atop \cmo}+{\cmo\atop \cpo}.$
\end{minipage}

\item[\rm(2C) and (2D)] follows from (2A), (2B) and Lemma \ref{duality}
\end{description}

\psfig{figure=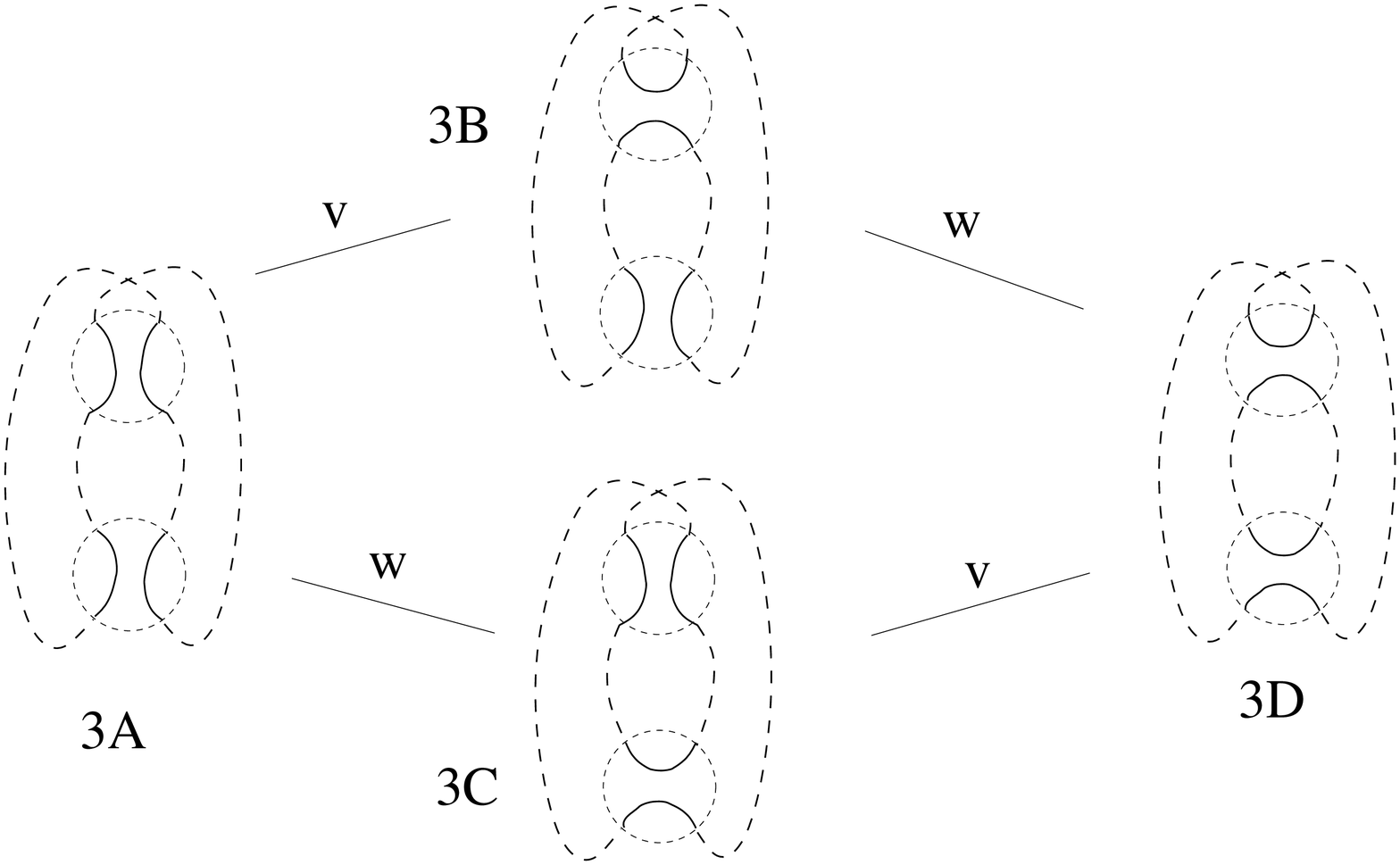,height=2in}

\begin{description}
\item[\rm(3A)] If 3D is NN, then both sides of eq. (\ref{dvdw}) vanish by
(\t). If 3D $\in$ TT then 3B, 3C are T and (\ref{dvdw}) follows by
(S). Analogously, if 3D is TN then 3B, 3C are N and (\ref{dvdw}) follows by
(S).
\item[\rm(3B)] By (NZ,\t) we assume that both 3B and 3C are trivial and
$S=$\cm. Hence 3D is either in TT or NN.
If 3A is nontrivial then $d_vd_w(\cm)=d_wd_v(\cm)=2\cp.$
Assume now that 3A is trivial. Since 3B, 3C are trivial as well,
$F=\RP^2.$ Since $d_wd_v(\cm)=0,$ $d_vd_w(\cm)=2\cp,$
(\ref{dvdw2}) holds.
\item[\rm(3C) and (3D)] follows from (3A), (3B) and Lemma \ref{duality}.
\end{description}

\centerline{\psfig{figure=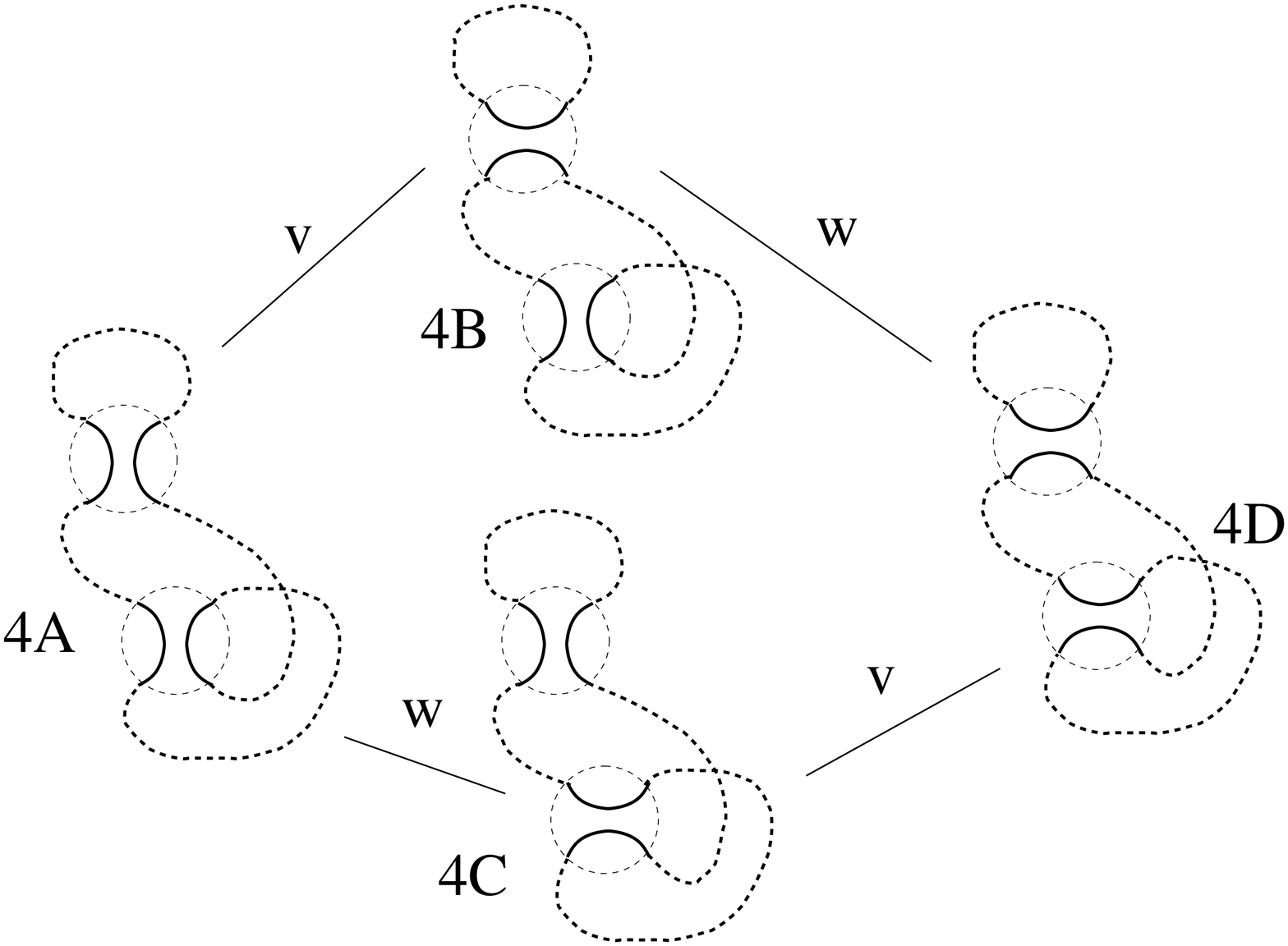,height=2in}}  

\begin{description}
\item[\rm(4A)] By (\t, NZ), (4A,4D)=(T,NT) or (N,TT). If 4A=T, then 4B=TT
or NN.  The following cases need to be considered:

$\diagfoura{T}{T}{T}{T}{N}{N}$ \begin{minipage}{4in}
In this case, $d_vd_w(S)=d_wd_v(S)$ is either ${\cp\atop \cpo}+
{\cp\atop \cmo}$

or $0$ depending on whether $S=\cm$ or \cp.
\end{minipage}

$\diagfoura{T}{N}{N}{N}{T}{N}$ \begin{minipage}{4in}
In this case, $d_vd_w(S)=d_wd_v(S)$ is either ${\cpo\atop \cp}+
{\cmo\atop \cp}$

or $0$ depending on whether $S=\cm$ or \cp.
\end{minipage}

$\diagfoura{N}{T}{N}{T}{T}{T}$ \begin{minipage}{4in}
In this case, $d_vd_w(S)=d_wd_v(S)={\cp\atop \cp}$ for
$S=\ceo.$
\end{minipage}

\item[\rm(4B)] By (\t,NZ), (4B,4C) is either (TT,N) and S=\cm\cm or 
(TN,T) and S=\cm\ceo.
In the first case, 4A=T and 4D=${T\atop N}.$
Hence $d_vd_w(\cm\cm)=d_wd_v(\cm\cm)=\cpo+\cmo.$
Assume now that 4B=TN, 4C=T. If the top circle in 4B is T then the top
circle in 4D is also T and since 4C is T, also the
bottom circle of 4D is T. Hence we have
$\diagfoura{N}{T}{N}{T}{T}{T}$ and in this case
$d_vd_w(S)=d_wd_v(S)$ is $\cp$ for S=\cm, and $0$ otherwise.

If the top circle in 4B is N then the top
circle in 4D is also N and since 4C is T, the
bottom circle of 4D is N. Hence we have
$\diagfoura{N}{N}{T}{N}{N}{T}$ and in this case
$d_vd_w(S)=d_wd_v(S)$ is $\cp$ for S=\cm, and $0$ otherwise.  

\item[\rm(4C) and (4D)] follows from (4A), (4B) and Lemma \ref{duality}.
\end{description}

\centerline{\psfig{figure=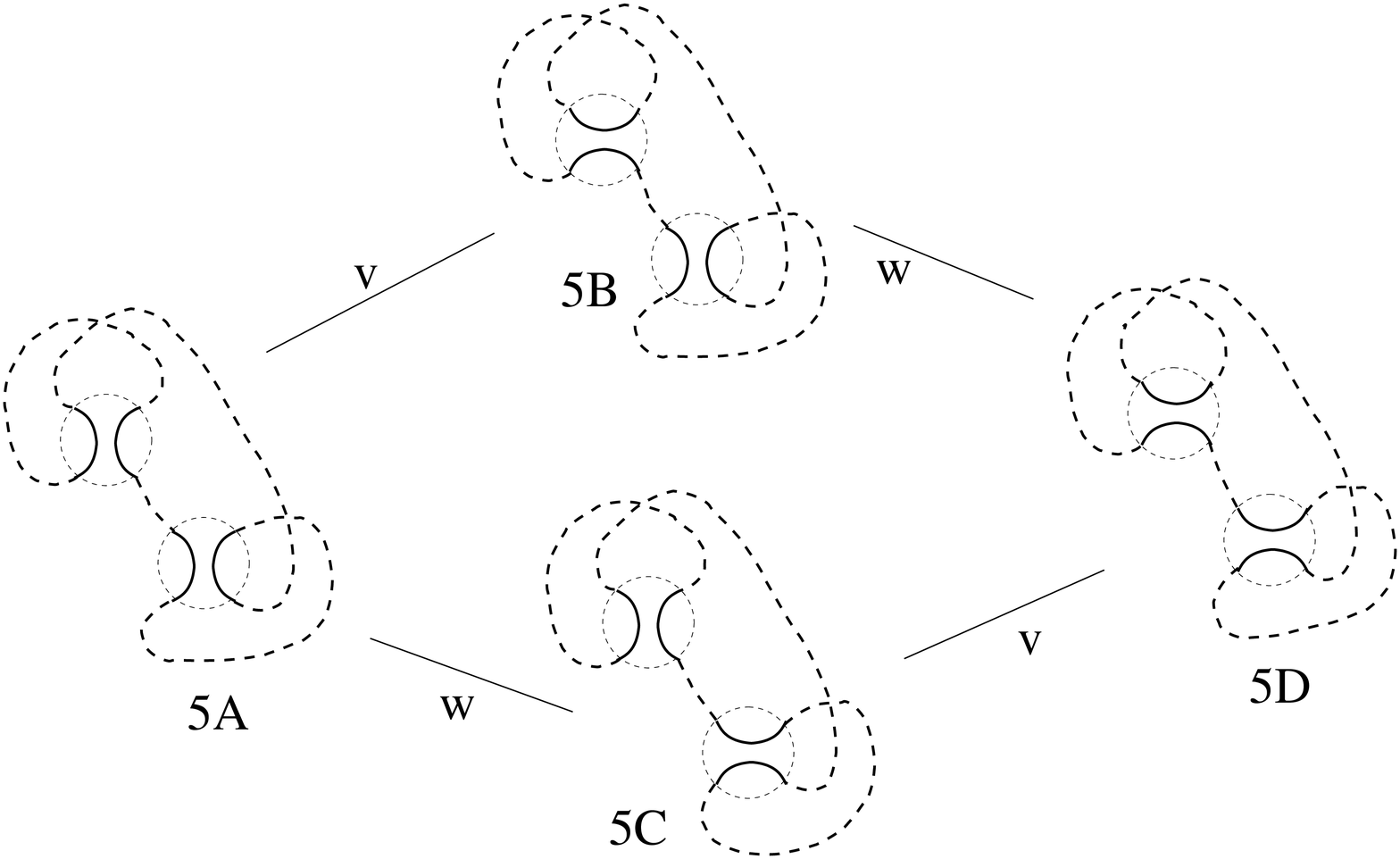,height=2in}}  

\begin{description}
\item[\rm(5A)] By (\t, NZ), we assume that both 5A and 5D are trivial.
Since for (5B,5C)=(T,T) or (N,N), (\ref{dvdw}) holds by (S),
assume that (5B,5C)=(T,N). Since (5A,5B)=(T,T), $F=\RP^2.$
In this case, $d_wd_v=0,$ and $d_vd_w(S)=0$ mod $2$ for all $S.$
\item[\rm(5B),(5C),(5D)] The proof is identical.
\end{description}

This completes the proof of Theorem \ref{d^2}.

%
\subsection{Reidemeister move I}
\label{r_I}
%

Let the link diagrams $D,D'$ differ by a negative kink only, 
$D=$\parbox{.8cm}{\psfig{figure=kinksmooth.eps, height=0.6cm}},$D'=$
\parbox{.8cm}{\psfig{figure=kink1.eps,height=0.6cm}} and let the map
$\rho_I:C_{ijs}(
\parbox{.8cm}{\psfig{figure=kinksmooth.eps,height=0.6cm}})\to
C_{i-1, j-3,s}(\parbox{.8cm}{\psfig{figure=kink1.eps,height=0.6cm}})$
be given by $\parbox{.8cm}{\psfig{figure=kinksmoothe.eps,height=.6cm}}\to 
\parbox{.6cm}{\psfig{figure=kink1nege.eps,height=.6cm}}$\hspace*{.1in} , where 
$\epsilon=\pm,\pm0.$ One easily verifies that $\rho_I$ is a chain map.
We are going to prove $\rho_I$ that induces an isomorphism on homology
groups by induction on the number of crossings of $D.$

Assume first that $D$ has no crossings. It this situation one may assume 
that $D$ is a single loop. Let $s=D\in \C(F)$ if $D$ is non-bounding
in $F,$ and let $s=0$ otherwise. Since $C_{**s'}(D)=C_{**s'}(D')=0,$ for
$s'\ne s,$ it is enough to consider chain groups $C_{**s}(D),$
$C_{**s}(D')$ only. Let $H_i(D)=\oplus_{j} H_{ijs}(D).$


By direct calculation, $H_i(D)=0$ for $i\ne 0,$ $H_i(D')=0$ for 
$i\ne -1,$ and $H_0(D)$ and $H_{-1}(D')$ have bases $\{\cpe,\cme\}$
and $\{\cpe\cm, \cme\cm\}$ respectively, where $\epsilon$ is either $1$ or $0$
depending on whether $D$ is contractible or not.
Since $\rho(\cpe)=\cpe\cm,$ $\rho(\cme)=\cme\cm,$ 
$\rho_I$ induces an isomorphism $\rho_{I*}:H_0(D)\to
H_{-1}(D'),$ preserving the grading with respect to $j$ and $s.$

Now assume that $D$ is a link diagram with $c>0$ crossings and 
that $\rho_{I*}$ is an isomorphism on homology groups for all 
diagrams with at most $c-1$ crossings. Let $D'$ be obtained from $D$
by adding a negative kink, as before. By inductive assumption,
the long exact sequence of homology groups of $D$ for any crossing
$p$ and the long exact sequence of homology groups of $D'$ for the same
$p$ are related by the maps $\rho_{I*}$ which, except for the central one,
are isomorphisms by the inductive assumption:
$$\begin{array}{ccccc}
\hspace*{-.1in} H_{i+2,j+2}(D_0) \stackrel{\partial}{\to} 
&\hspace*{-.1in} H_{i+2,j+4}(D_\infty)
\stackrel{\al_*}{\to} & \hspace*{-.1in} H_{i+1,j+3}(D)
\stackrel{\be_*}{\to}& \hspace*{-.1in}
H_{i,j+2}(D_0) \stackrel{\partial}{\to} & \hspace*{-.1in}
H_{i,j+4}(D_\infty)\\
\downarrow \rho_I  & \downarrow \rho_I & \downarrow \rho_I &
\downarrow \rho_I  & \downarrow \rho_I\\
\hspace*{-.1in} H_{i+2,j+2}(D_0') \stackrel{\partial}{\to} & \hspace*{-.1in} 
H_{i+2,j+4}(D_\infty')
\stackrel{\al_*}{\to} & \hspace*{-.1in} H_{i+1,j+3}(D') \stackrel{\be_*}{\to}& 
\hspace*{-.1in} H_{i,j+2}(D_0') \stackrel{\partial}{\to} &
\hspace*{-.1in} H_{i,j+4}(D_\infty').
\end{array}$$
(For brevity, we skip the index ``s''.)

Since for any state $S$ of $D,$ $\rho_I(S)$ differs from $S$ by an
additional \cm only, and that \cm does not pass trough any crossings
of $D,$ the above diagram commutes.
Therefore, by Five lemma, $\rho_I: H_{***}(D)\to H_{***}(D')$ is an
isomorphism and the proof of Theorem \ref{main_thm}(2) 
is completed.

%
\subsection{Invariance under Reidemeister move II}
\label{r_II}
%

Consider two link diagrams \pic{L0}, 
\parbox{.8cm}{\psfig{figure=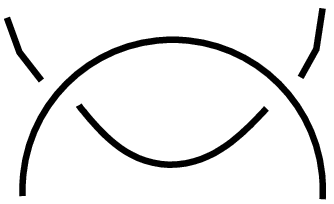,height=.5cm}} related to each other 
by the second Reidemeister move and denoted by $D$ and $D'$
respectively.
These two diagrams differ by two vertices which will be denoted by $v$
and $w,$ \parbox{1.2cm}{\psfig{figure=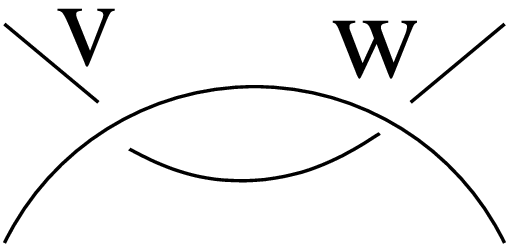,height=.5cm}}.
We assume that the vertices of $D'$
are ordered such that
$v$ is the first and that $w$ is the second
vertex and that the vertices of $D$ are ordered by
an ordering inherited from the ordering of vertices of $D'.$

Let $f: C_{ijs}(\pic{L0})\to 
C_{ijs}(\parbox{.8cm}{\psfig{figure=r2.eps,height=.5cm}})$ and
$g: C_{ijs}(\pic{Linfty})\to 
C_{i,j-2,s}(\parbox{.8cm}{\psfig{figure=r2.eps,height=.5cm}})$
be defined as follows:
$$\pic{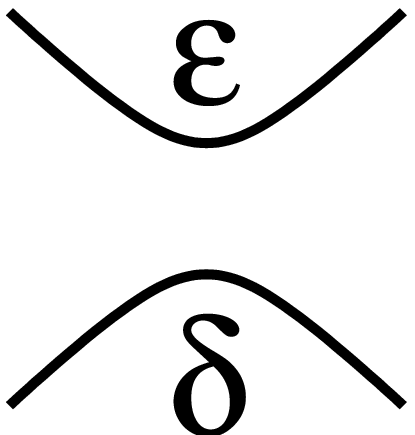}\ \stackrel{f}{\to}\ 
\parbox{1.1cm}{\psfig{figure=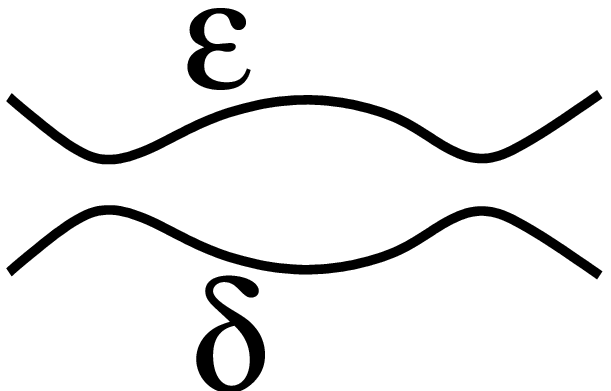,height=.7cm}},\quad
\pic{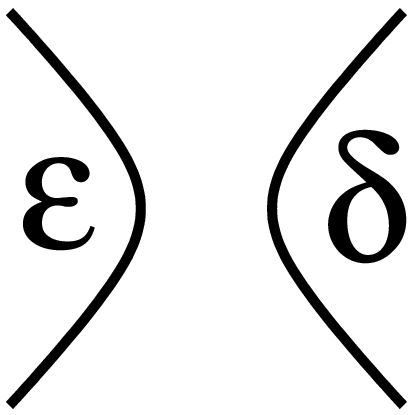}\ \stackrel{g}{\to}\ 
\parbox{1.2cm}{\psfig{figure=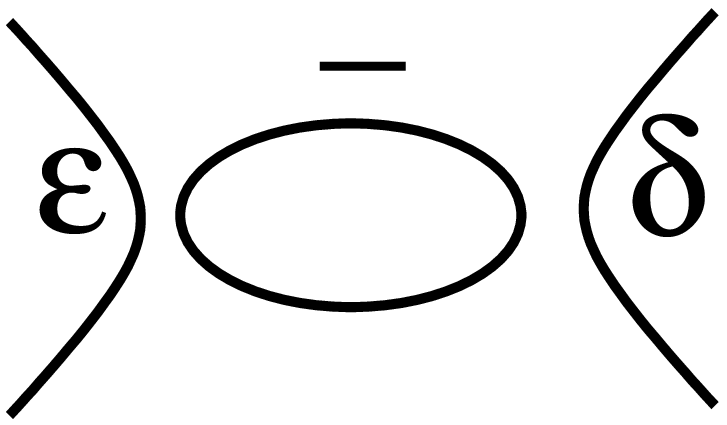,height=.6cm}}.$$
Let $\rho_{II}: C_{ijs}(\pic{L0})\to 
C_{ijs}(\parbox{.8cm}{\psfig{figure=r2.eps,height=.5cm}})$
be given by $$\rho_{II}(S)=f(S)+g(\gamma(S)).$$
(Recall that $\gamma: C_{ijs}(\pic{L0}) \to
C_{i,j+2,s}(\pic{Linfty})$ was defined by
(\ref{gamma}) in Section \ref{s_viro}.) 

\begin{theorem}\label{r2_thm}
{\rm(1)}\qua $\rho_{II}$ is a chain map.

{\rm(2)}\qua $\rho_{II*}:H_{ijs}(\pic{L0})\to 
H_{ijs}(\parbox{.8cm}{\psfig{figure=r2.eps,height=.5cm}})$
is an isomorphism.
\end{theorem}

We start the proof with the following lemma whose proof is left to the reader.

\begin{lemma}\label{r2_lemma}
{\rm(1)}\qua For any enhanced state $S\in C_{i,j+2,s}(\pic{Linfty}),$
$dg(S)=gd(S)+(-1)^{m(S)+1}\imath(S),$ where $\imath(S)\in
  C_{i-2,j,s}(\parbox{.8cm}{\psfig{figure=r2.eps,height=.5cm}})$
is the enhanced state of the form \pic{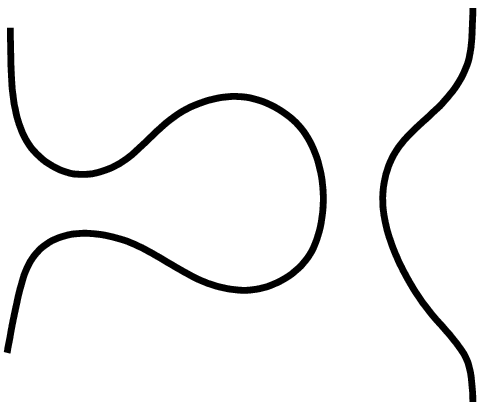}
obtained by ``deforming'' the diagram of $S.$
In other words, the labels of $S$ and $\imath(S)$ coincide under the obvious
$1-1$ correspondence of circles of $S$ and $\imath(S).$

{\rm(2)}\qua For any enhanced state $S\in C_{ijs}(\pic{L0}),$
$df(S)=fd(S)+ (-1)^{m(S)}\sum S',$ where the sum is taken over all 
enhanced states
$S'\in C_{i-2,j,s}(\parbox{.8cm}{\psfig{figure=r2.eps,height=.5cm}})$
of the form \pic{r2smooth2}
and such that the labels of all circles of $S$ and of $S'$ coincide,
except possibly the circles which pass through \pic{L0} in $S$ and
through \pic{r2smooth2} in S'.
\end{lemma}

By Proposition \ref{gamma_prop}(4), the sum of all states $S'$ 
defined above equals to $\imath(\gamma(S)).$ Hence
\begin{equation}\label{eq_f}
df(S)=fd(S)+(-1)^{m(S)}\imath \gamma(S).
\end{equation}
If $\gamma(S)=\sum_i S_i$ then $m(S_i)=m(S)$ for all $i,$ and therefore
by Lemma \ref{r2_lemma}(1)
$$dg\gamma(S)=gd\gamma(S)+(-1)^{m(S)+1}\imath\gamma(S)$$
and, by Proposition \ref{gamma_prop}(3),
\begin{equation}\label{eq_g}
dg\gamma(S)=
g\gamma d(S)+(-1)^{m(S)+1}\imath\gamma(S).
\end{equation}
By adding (\ref{eq_f}) and (\ref{eq_g}) we get the statement of
Theorem \ref{r2_thm}(1).

The remaining part of this section is devoted to the proof of Theorem
\ref{r2_thm}(2).

For any enhanced state $S,$ let $m(S)$ denote the number of negative 
crossing markers in $S$ and let $\eta(S)=(-1)^{m(S)}S.$ 

\begin{lemma} $d\eta=-\eta d$
\end{lemma}

\begin{proof} It is enough to prove that $\hat d_p\eta(S)=-\eta\hat
  d_p(S)$ for each crossing $p$ in every state $S.$ 
If the marker at $p$ is negative then both sides of the above equation
vanish. Otherwise, $d_p(S)=\sum S'$ and $m(S')=m(S)+1$ for every
summand $S'$ in $d_p(S).$ 
\end{proof}


\begin{lemma}\label{1st_square}
The following diagram commutes:
$$\begin{array}{ccc}
C_{ijs}(\picll{D--})
 & \stackrel{\alpha}{\to} &
C_{i-1,j-1,s}(\picll{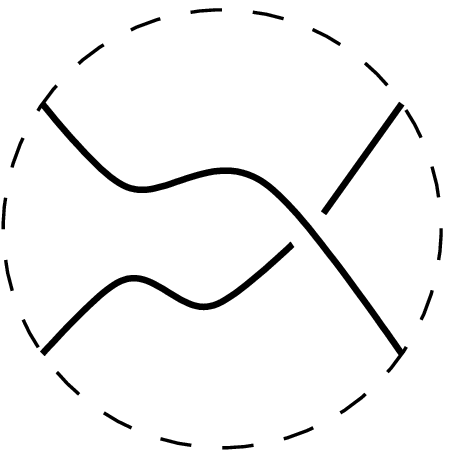})\\
\downarrow \rho_I &  & \downarrow \eta \\
C_{i-1,j-3,s}(\picll{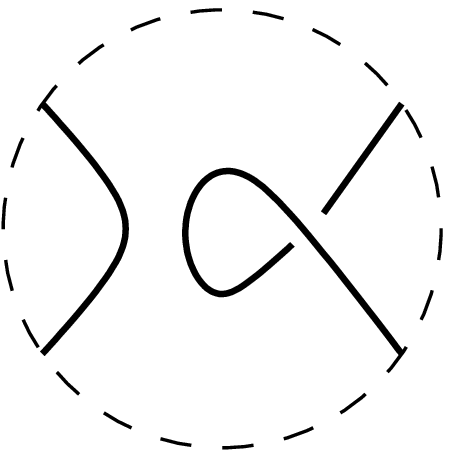})& \stackrel{\hat \gamma}{\to} &
C_{i-1,j-1,s}(\picll{D-0})
\end{array}$$
\end{lemma}

\proof By definition of $\alpha$ given in Section \ref{s_viro},
$\alpha(S)=\alpha_0(S)$ for any enhanced state $S\in
  C_{ijs}(\parbox{.8cm}{\psfig{figure=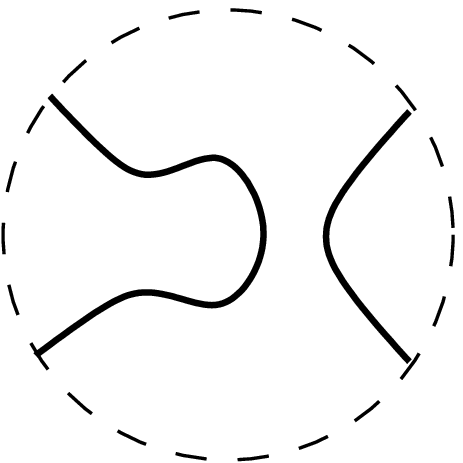,height=.6cm}}).$
Furthermore
$\rho_I(\parbox{.8cm}{\psfig{figure=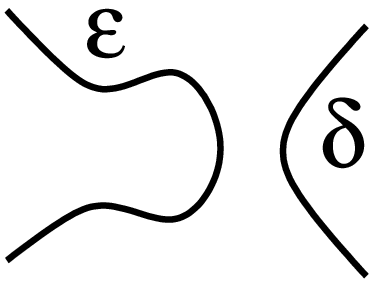,height=.6cm}})=
\parbox{1.2cm}{\psfig{figure=Linftyc.eps,height=.6cm}}$
and $\gamma(\parbox{1.2cm}{\psfig{figure=Linftyc.eps,height=.6cm}})=
\parbox{.8cm}{\psfig{figure=D--l.eps,height=.6cm}},$ 
by Proposition \ref{gamma_prop}(4).
Hence $\gamma \rho_I=\alpha,$ and 
$$\hat \gamma \rho_I(S)=(-1)^{m(S)+1}\gamma
\rho_I(S)=(-1)^{m(S)+1}\alpha(S) =\eta \alpha(S).\eqno{\qed}$$

\begin{lemma}
The following diagram commutes:
$$\begin{array}{ccc}
H_{i+1,j+1,s}(\picll{D-0}) & \stackrel{\beta_*}{\to} &
H_{ijs}(\picll{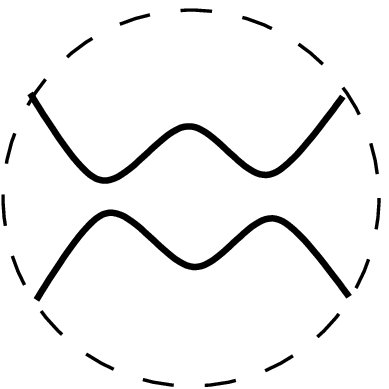})\\
  & \searrow \alpha_* & \downarrow \rho_{II*}\\
 & & H_{ijs}
(\picll{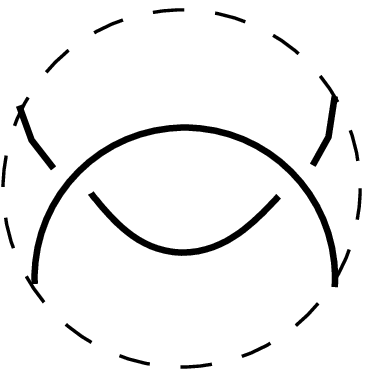})
\end{array}$$
\end{lemma}

\proof Any $x\in C_{i+1,j+1,s}(\pic{D-0})$ decomposes into
$x^+ +x^-$ where $x^+$ (respectively: $x^-$) is a linear combination 
of states with positive (respectively: negative) marker at the
crossing $w.$ Assume that $x\in Z_{i+1,j+1,s}(\pic{D-0}).$
Since $\beta(x^-)=0,$ we need to show that
\begin{equation}\label{eq_boundary}
(\rho_{II}\beta-\alpha)(x^+
+x^-)=f\beta(x^+)+g\gamma\beta(x^+)-\alpha(x^+)-\alpha(x^-)
\end{equation}
is a boundary cycle.
Since the diagram 
$$\begin{array}{ccc}
C_{i+1,j+1,s}(\picll{D-0}) & \stackrel{\beta}{\to} &
C_{ijs}(\picll{r2smoothc})\\
  & \searrow \alpha_0=\alpha & \downarrow f\\
 & & C_{ijs}
(\picll{r2c})
\end{array}$$
commutes for $x^+$, the right hand side of (\ref{eq_boundary}) equals
$g\gamma\beta(x^+)-\alpha(x^-).$
Since $$\gamma\beta=\hat \gamma \eta \beta=
\bar \alpha \hat d_w \bar \beta \eta \beta = 
- \bar \alpha \eta \hat d_w \bar \beta \beta=
- \bar \alpha \eta \hat d_w,$$ 
\begin{equation}\label{eq_comm}
(\rho_{II}\beta-\alpha)(x)=
-g\bar \alpha \eta \hat d_w(x^+)-\alpha(x^-).
\end{equation}
Since $x=x^++x^-\in Ker\, d,$ $$\sum_c \hat d_c(x^+)+
\sum_c \hat d_c(x^-)=0,$$ and by taking the states with the
negative marker at $w$ into account only, we get
$\hat d_w(x^+)+ \sum_c \hat d_c(x^-)=0.$
By substituting $-\sum_c \hat d_c(x^-)$ for $d_w(x^+)$
in (\ref{eq_comm}) we get
$$(\rho_{II}\beta-\alpha)(x)=
g\bar \alpha \eta\left( \sum_c \hat d_c(x^-)\right) -\alpha(x^-),$$
and, since $g$ and $\bar \alpha$ commute with $\hat d_c$ for $c\ne
v,w,$ (cf. (\ref{alphabeta_com})), and $\eta \hat d_c=- \hat d_c \eta,$
$$(\rho_{II}\beta-\alpha)(x)=
-\sum_{c\ne v,w} \hat d_c(g \bar\alpha\eta(x^-)) -\alpha(x^-).$$
If $x^-=\sum_i a_i S_i$ and $S_i=
\parbox{1.1cm}{\psfig{figure=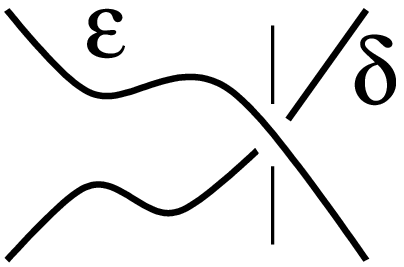,height=.6cm}}$
then $$g\bar \alpha \eta(S_i)= 
\parbox{1.2cm}{\psfig{figure=Linftyc.eps,height=.6cm}}\cdot
(-1)^{m(S_i)}\in C_{ijs}(\parbox{.8cm}{\psfig{figure=r2.eps,height=.5cm}}).$$
Therefore $$d(g\bar \alpha \eta(S_i))=\sum_{c\ne v,w} 
\hat d_c(g \bar\alpha\eta(S_i)) + 
\parbox{.8cm}{\psfig{figure=D--l.eps,height=.6cm}}\cdot (-1)^{m(S_i)}
(-1)^{m(S_i)}.$$
Hence
$$(\rho_{II}\beta-\alpha)(x)=-d(g\bar \alpha \eta(x^-)).\eqno{\qed}$$

\begin{corollary}\label{2nd_square}
The following diagram commutes:
$$\begin{array}{ccc}
H_{i+1,j+1,s}(\picll{D-0}) & \stackrel{\beta_*}{\to} &
H_{ijs}(\picll{r2smoothc})\\
\downarrow \eta_* &  & \downarrow \rho_{II*}\eta_*\\
H_{i+1,j+1,s}(\picll{D-0})& \stackrel{\alpha_*}{\to}& H_{ijs}
(\picll{r2c})
\end{array}$$
\end{corollary}

\begin{proof}
Since $\eta$ skew-commutes with $\alpha,$ $f,$ and $g\gamma,$
we have 
$$\alpha_*\eta_*=-\eta_*\alpha_*=-\eta_*\rho_{II*}\beta_*=
\rho_{II*}\eta_*\beta_*,$$
by the previous lemma.
\end{proof}

\begin{lemma}\label{3rd_square}
The following diagram commutes:
$$\begin{array}{ccc}
C_{ijs}(\picll{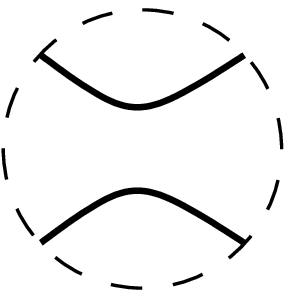}) 
& \stackrel{\hat\gamma}{\to} &
C_{i,j+2,s}(\picll{D--})\\
\downarrow \rho_{II}\eta&  & \downarrow \rho_I \\
C_{ijs}(\picll{r2c})&
\stackrel{\beta}{\to} & C_{i-1,j-1,s}(\picll{D+0}) 
\end{array}$$
\end{lemma}

\begin{proof} Since $\beta f(S)=0$ for $S\in C_{ijs}(\pic{r2Dp}),$
$\beta \rho_{II}\eta(S)=(-1)^{m(S)}\beta g \gamma(S).$ 
Hence we need to prove that 
$$(-1)^{m(S)}\beta g \gamma(S)=(-1)^{m(S)}\rho_I\gamma(S).$$
Now the statement follows
from the following commutative diagram:
$$\begin{array}{ccc}
  & & C_{i,j+2,s}(\picll{D--})\\
 & \swarrow g  & \downarrow \rho_I \\
C_{ijs}(\picll{r2c})&
\stackrel{\beta}{\to} & C_{i-1,j-1,s}(\picll{D+0}) 
\end{array}$$

\vspace{-25pt}
\end{proof}
\medskip

Consider the following two Viro's exact sequences
$$\begin{array}{c@{}c@{}c@{}c@{}c}
 H_{i+2,j+2}(\pic{D--}) \stackrel{\al_*}{\to} & H_{i+1, j+1}(\pic{D-0}) 
  \stackrel{\be_*}{\to} &  H_{ij}(D) \stackrel{\partial}{\to} &
  H_{i,j+2}(\pic{D--})  \stackrel{\al_*}{\to}  &  H_{i-1,j+1}(\pic{D-0}) \\
 \downarrow \rho_{I*} & \downarrow \eta_* & \downarrow \rho_{II*}\eta_* &
 \downarrow \rho_{I*} & \downarrow \eta_*  \\
 H_{i+1, j-1}(\pic{D+0}) \stackrel{\partial}{\to} & H_{i+1,
 j+1}(\pic{D-0}) \stackrel{\al_*}{\to} & H_{ij}(D')
 \stackrel{\be_*}{\to} & H_{i-1, j-1}(\pic{D+0}) \stackrel{\partial}{\to} &
 H_{i-1,j+1}(\pic{D-0})
\end{array}  $$
(The third index, $s,$ is suppressed to shorten notation.)

By Lemma \ref{1st_square}, Corollary \ref{2nd_square}, and Lemma 
\ref{3rd_square}, the above diagram commutes.
Hence, by Five Lemma, $\rho_{II*}\eta_*$ is an isomorphism.
Since $\eta\eta=id,$ $\eta_*$ is an isomorphism, and hence
$\rho_{II*}$ is an isomorphism as well.
This completes the proof of Theorem \ref{r2_thm}(2).

%
\subsection{Invariance under Reidemeister move III}
\label{r_III}
%

In this section we prove the invariance of Khovanov homology under
third Reidemeister moves. The proof is identical for homologies based on
trivial and separating circles.
Let $D, D'$ be link diagrams related to each other by a third
Reidemeister move shown below:
$$\parbox{1.5cm}{\psfig{figure=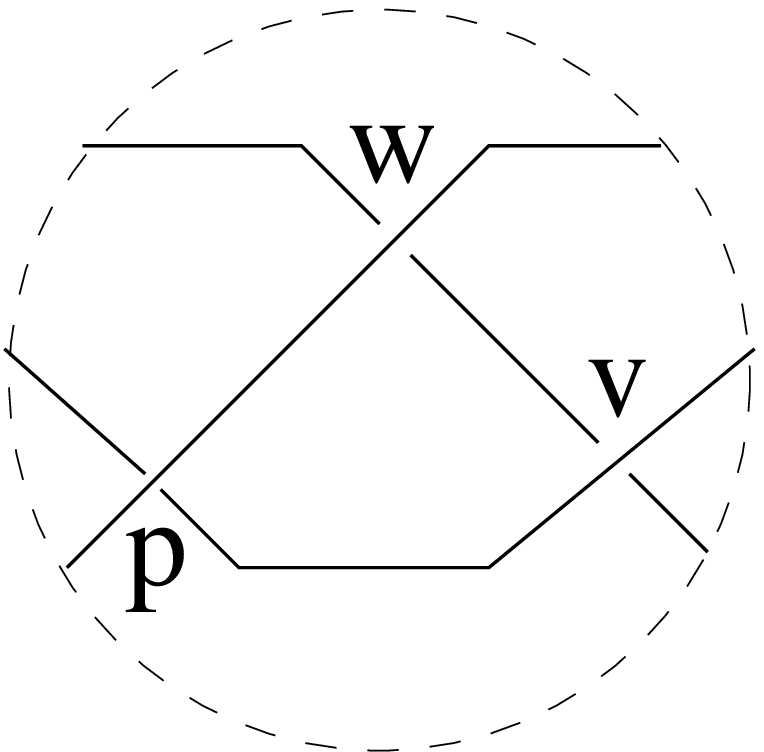,height=1.5cm}}\subset D\quad 
{\rm and}\quad \parbox{1.5cm}{\psfig{figure=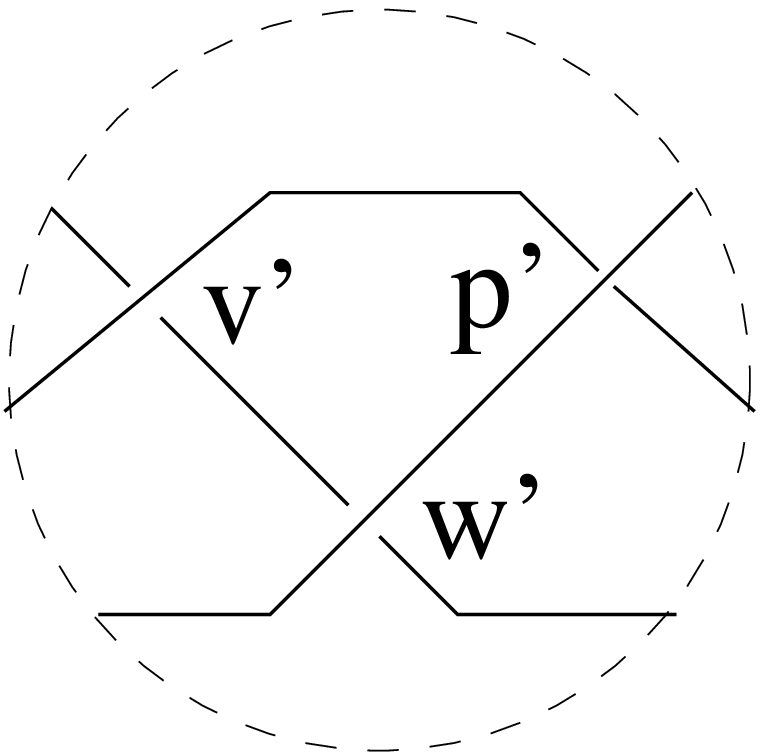,height=1.5cm}} 
\subset D'.$$
Assume that the vertices of $D$ and $D'$ are ordered such that 
$p,v,w$ (respectively: $p',v',w'$) are the first, the second, and the
third among vertices of $D$ (respectively: of $D'$) and that the orderings
of the remaining vertices of $D$ and $D'$ coincide.
We denote by $D_\pm$ (respectively: by $D_\pm'$) the diagrams obtained from $D$ 
(respectively: $D'$) by smoothing the crossing $p$ (respectively:
$p'$) with $\pm$ marker. 
$$ \begin{array}{ccc}
 D_- & \cong& D'_- \\
  \cup && \cup \\
 \parbox{1.5cm}{\psfig{figure=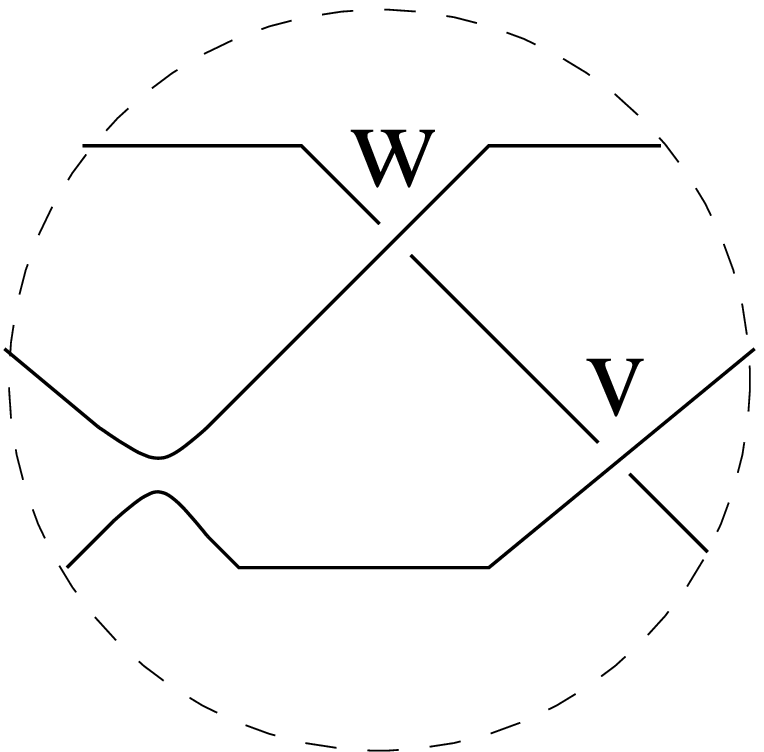,height=1.5cm}}\quad & \cong &\quad
    \parbox{1.5cm}{\psfig{figure=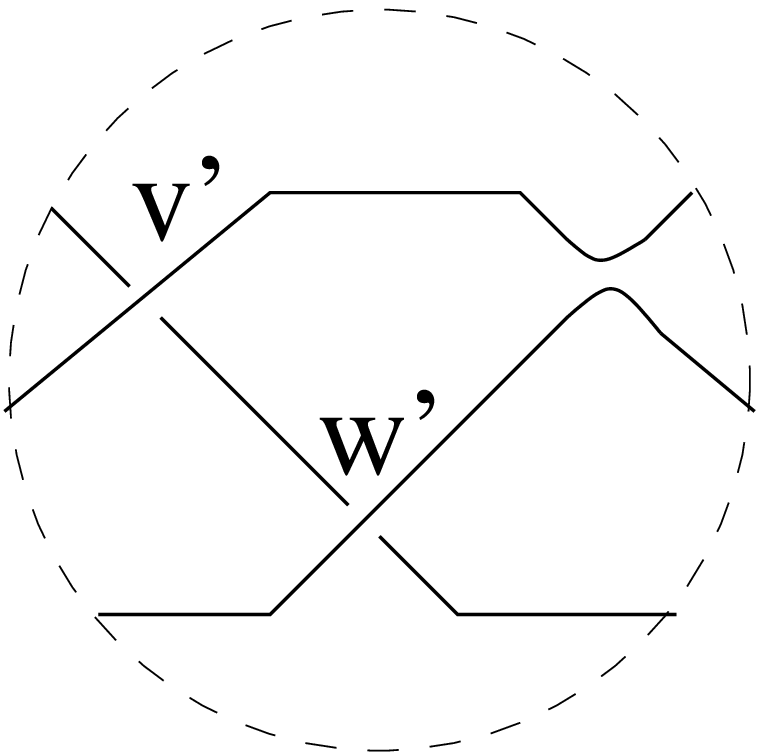,height=1.5cm}}
  \end{array} $$ $$
  \begin{array}{ccccccc}
  D_+ &  \stackrel{\rm R2}{\Longleftrightarrow}  & D_{++-} & \cong & D'_{++-} 
 & \stackrel{\rm R2}{\Longleftrightarrow} & D'_+ \\
\cup & & \cup & & \cup & & \cup \\
\parbox{1.5cm}{\psfig{figure=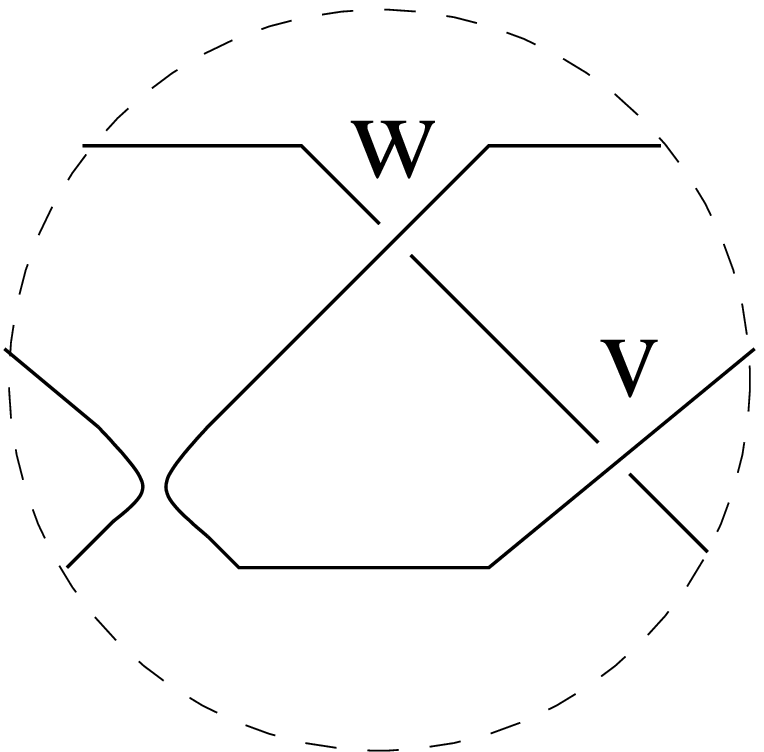,height=1.5cm}}\quad & 
\stackrel{\rm R2}{\Longleftrightarrow}\quad &  
\parbox{1.5cm}{\psfig{figure=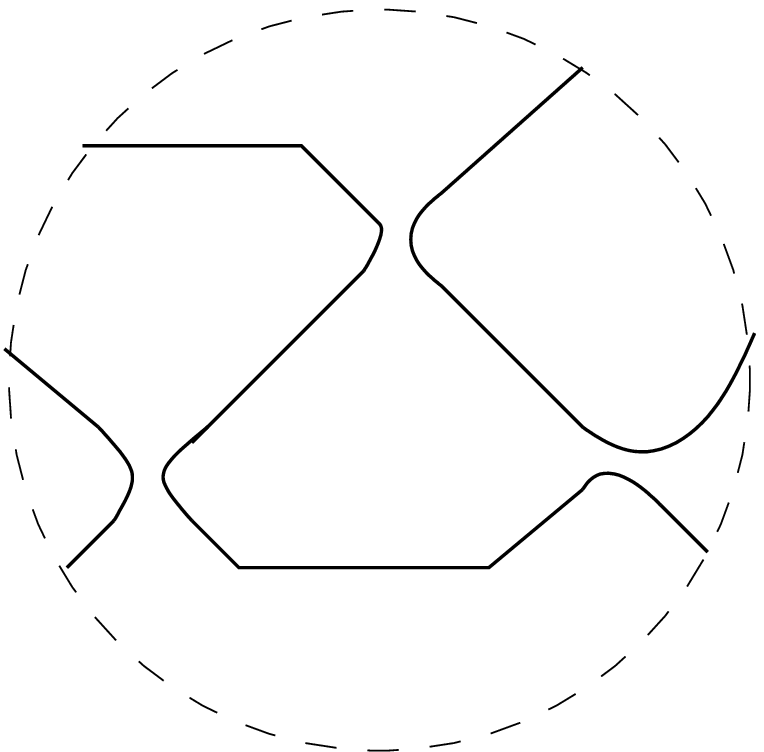,height=1.5cm}}\quad & \cong\quad &
\parbox{1.5cm}{\psfig{figure=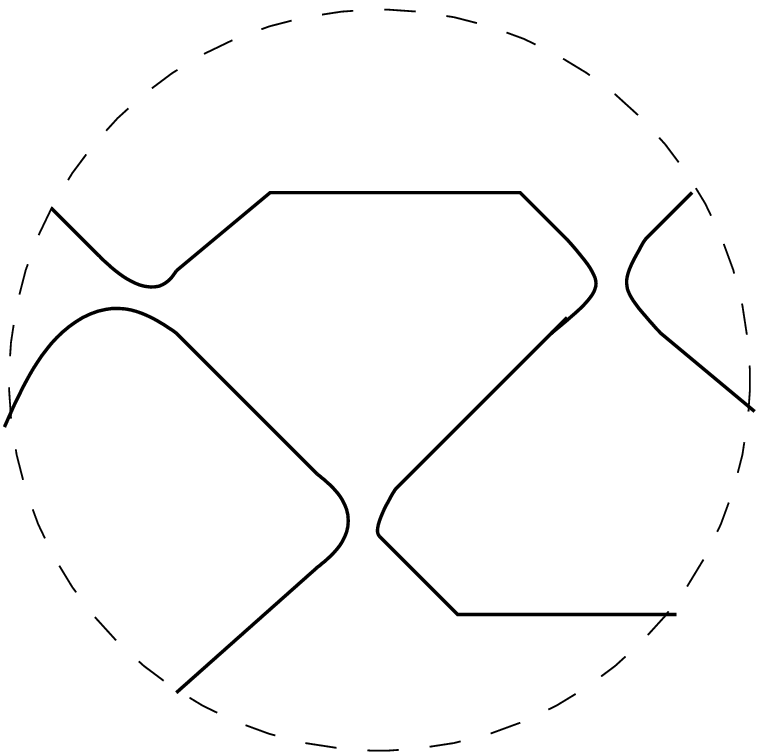,height=1.5cm}}\quad & 
   \stackrel{\rm R2}{\Longleftrightarrow}\quad&
    \parbox{1.5cm}{\psfig{figure=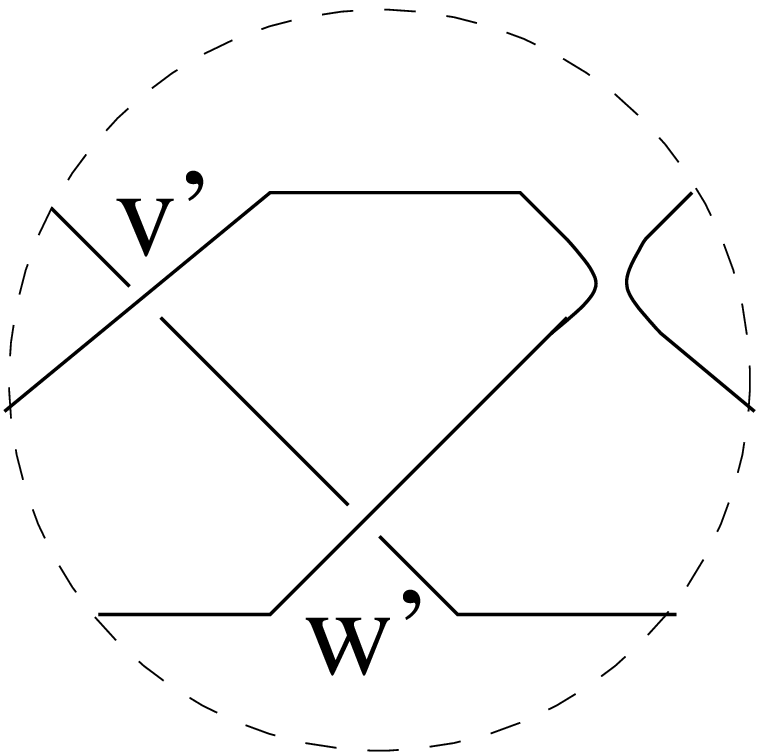,height=1.5cm}},
\end{array} $$
Let $\rho_{II}: C_{***}(D_{++-}) \to C_{***}(D_+)$, and 
$\rho_{II}': C_{***}(D'_{++-}) \to C_{***}(D'_+)$ be the isomorphisms 
defined in Section \ref{r_II}.

Since the diagrams $D_-$ and $D_-'$ are identical, except the crossing
ordering
change between vertices $v$ and $w,$ their homology groups are isomorphic
via the isomorphism $f: H_{ijs}(D_-)\to H_{ijs}(D_-')$ defined in the
Section \ref{sec_independence}.
\comment{
If $S_{\ve_v,\ve_w}$ is an enhanced state of $C(D_-)$ with 
$\ve_v, \ve_w$ markers at $v$ and $w$ and
$S_{\ve_w,\ve_v}'$ denotes the corresponding enhanced state of $C(D_-')$ with 
$\ve_w, \ve_v$ markers at $w'$ and $v',$ respectively, then
$$f(S_{\ve_v,\ve_w})=\begin{cases} -S_{--}' & \text{if $\ve_v=\ve_w=-1$} \\
S_{\ve_w,\ve_v}' & \text{otherwise.}\\ \end{cases}$$
} 

Consider Viro's exact sequence:
$$0\to C_{i+1,j+1,s}(D_-)\stackrel{\alpha}{\to}
C_{ijs}(D)\stackrel{\beta}{\to} C_{i-1,j-1,s}(D_+)\to 0.$$
Let $C_{ijs}'(D_+)=\rho_{II}(C_{ijs}(D_{++-}))$ and
let $C_{ijs}'(D)=\beta^{-1}(C_{ijs}'(D_+)).$ 
\comment{ It is easy to prove that
$C_{***}'(D_+)$ is a sub-chain complex of $(C_{***}(D_+),d)$ and
$C_{***}'(D)$ is a sub-chain complex of $(C_{***}(D),d).$
} 

\begin{proposition}\label{embedding}
The embeddings

{\rm(1)}\quad $C_{***}'(D_+)\hookrightarrow C_{***}(D_+)$

{\rm(2)}\quad $C_{***}'(D)\hookrightarrow C_{***}(D)$

induce isomorphisms on homology groups.
\end{proposition}

\begin{proof}
(1) follows from Theorem \ref{r2_thm}(2).

(2)\qua Since $\alpha(C_{i+1,j+1,s}(D_-))\subset C_{ijs}'(D),$
Viro's exact sequence restricts to
$$0\to C_{i+1,j+1,s}(D_-)\stackrel{\alpha}{\to}
C_{ijs}'(D)\stackrel{\beta}{\to} C_{i-1,j-1,s}'(D_+)\to 0.$$
Therefore we have the following commuting diagram in which the
vertical maps are induced by the natural embeddings:
$$\begin{array}{c@{}c@{}c@{}c@{}c}
H_{i+1, j-1}'(D_+) \stackrel{\partial}{\to} & H_{i+1, j+1}(D_{-})  
\stackrel{\al_{*}}{\to} & H_{ij}'(D) \stackrel{\be_{*}}{\to} &
H_{i-1, j-1}'(D_+) \stackrel{\partial}{\to} & H_{i-1, j+1}(D_{-}) \\
\downarrow & \downarrow & \downarrow & \downarrow & \downarrow \\
H_{i+1, j-1}(D_+) \stackrel{\partial}{\to} & H_{i+1,j+1}(D_-) 
\stackrel{\al_*}{\to} & H_{ij}(D)
\stackrel{\be'_*}{\to} & H_{i-1,j-1}(D_+) \stackrel{\partial}{\to} &
 H_{i-1,j+1}(D_-).
\end{array}$$
 Since the vertical morphisms $H_{***}(D_-)\to
H_{***}(D_-)$ are identity maps, the statement follows by Five Lemma.
\end{proof}

Similarly, we define $C_{***}'(D_+')\subset C_{***}(D_+')$
and $C_{***}'(D')\subset C_{***}(D')$ and prove that 
these embeddings induce isomorphisms on homology groups.

Since $\rho_{II}: C_{ijs}(D_{++-})\to C_{ijs}'(D_+)$ is an isomorphism,
the map $\rho=\rho_{II}' \rho_{II}^{-1},$ 
$\rho: C_{***}'(D_+)\to C_{***}'(D_+'),$ is well defined.

\begin{proposition}\label{IIIcommutes}
 The following diagram commutes:
$$\begin{array}{ccccc}
0 \to & C_{i+1,j+1,s}(D_-) 
\stackrel{\al}{\to} & C_{ijs}'(D) \stackrel{\be}{\to} & 
C_{i-1,j-1,s}'(D_+) \to & 0\\
 & \downarrow f & \downarrow \rho_{III} & 
\downarrow \rho &   \\
0 \to & C_{i+1,j+1,s}(D_-') \stackrel{\al}{\to} & C_{ijs}'(D') 
\stackrel{\be}{\to} & C_{i-1,j-1,s}'(D_+') \to & 0,
\end{array} $$
where $\rho_{III}=\bar \beta\rho\be + \al f \bar\al$ and $\bar\alpha,
\bar\beta$ are the maps defined in Section \ref{s_viro}.

\end{proposition}

\begin{proof}
We have $\rho_{III}\alpha=(\bar \beta\rho\be + \al f \bar\al)\alpha=
\alpha f$ and $\beta\rho_{III}=\beta(\bar\beta\rho\be - \al f \bar\al)=\rho\be.$
\end{proof}

\begin{lemma}\label{partialcommutes}
The following diagram commutes:
$$\begin{array}{ccc}
C_{i+1,j-1,s}'(D_+) & \stackrel{\hat\gamma}{\to} & C_{i+1,j+1,s}(D_-)\\
 \downarrow \rho &  &\downarrow f\\
C_{i+1, j-1,s}'(D_+') & \stackrel{\hat\gamma}{\to} & C_{i+1,j+1,s}(D_-')
\end{array} $$
\end{lemma}

\begin{proof}
Since $\rho=\rho_{II}'\rho_{II}^{-1},$ we need to show that
\begin{equation}\label{first_sq}
f\hat\gamma \rho_{II}(S)=\hat\gamma \rho_{II}'(S)
\end{equation} for
$S\in C_{***}(D_{++-}).$
If $S=$ \parbox{1.5cm}{\psfig{figure=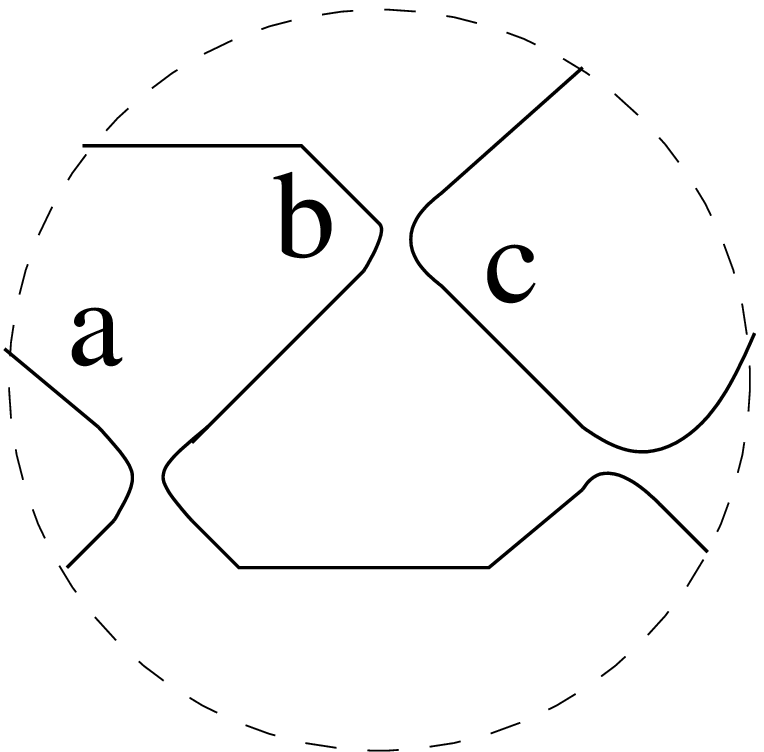,height=1.5cm}}
$\in C_{ijs}(D_{++-})$ then 
$$\rho_{II}(S)=\parbox{1.5cm}{\psfig{figure=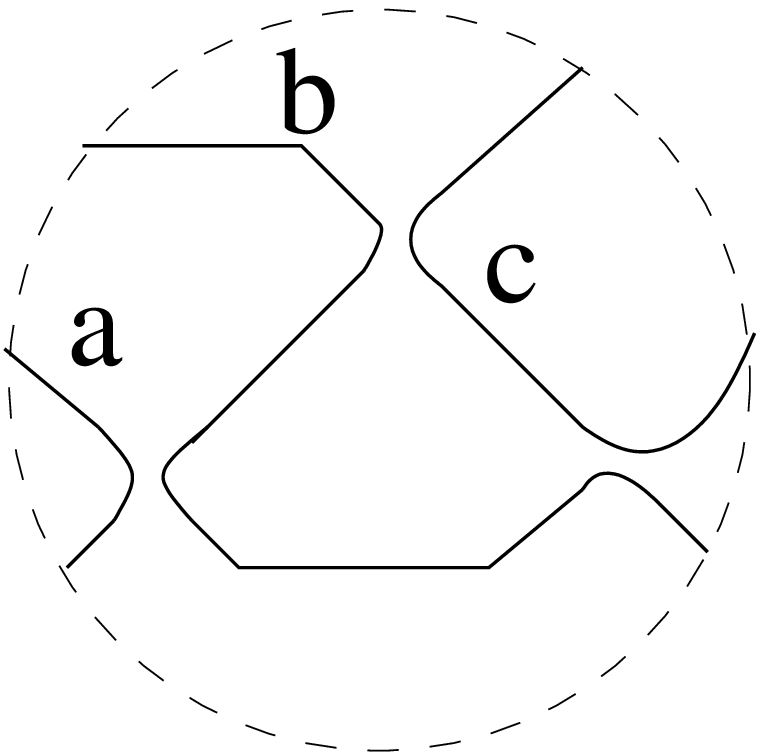,height=1.5cm}}
+ \sum \parbox{1.5cm}{\psfig{figure=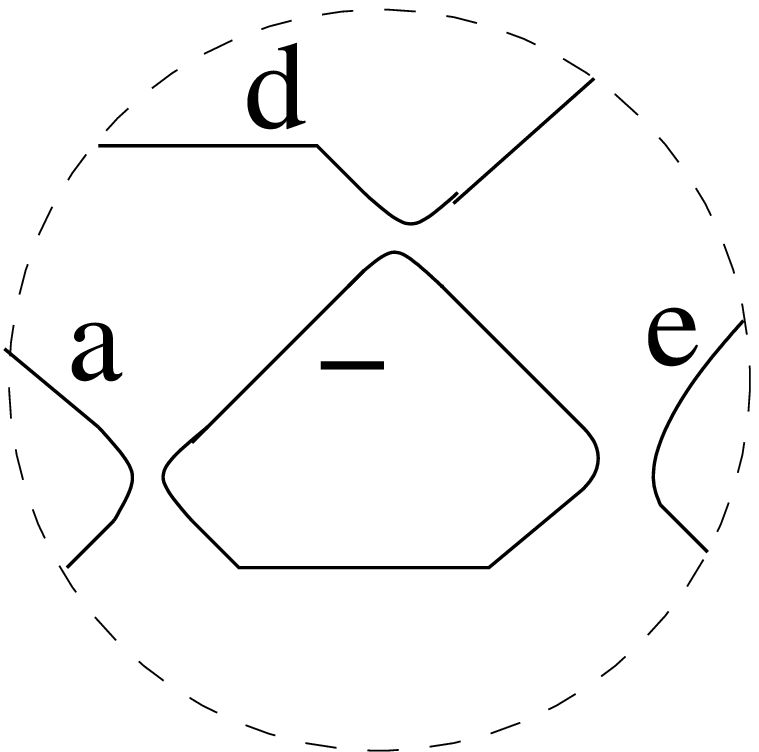,height=1.5cm}},$$
where the sum is over all possible labels $d,e$ for which the above
states are in $C_{ijs}(D_+).$
Therefore, by (\ref{gamma}) and Proposition \ref{gamma_prop}(4),
\begin{equation}\label{parrho}
\hat\gamma\rho_{II}(S)=(-1)^{m(S)+1}\sum 
\parbox{1.5cm}{\psfig{figure=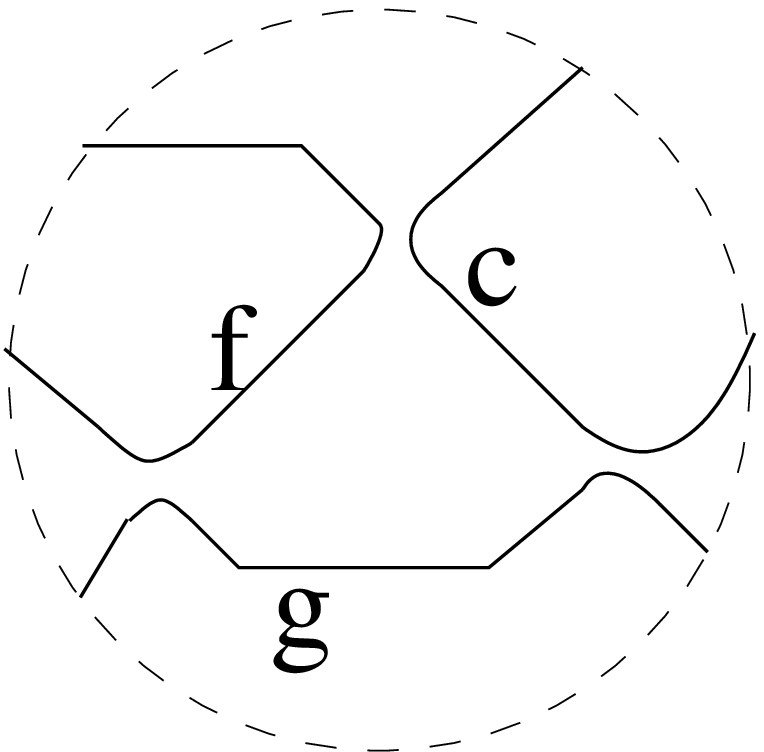,height=1.5cm}}
+ (-1)^{m(S)+1}\sum
\parbox{1.5cm}{\psfig{figure=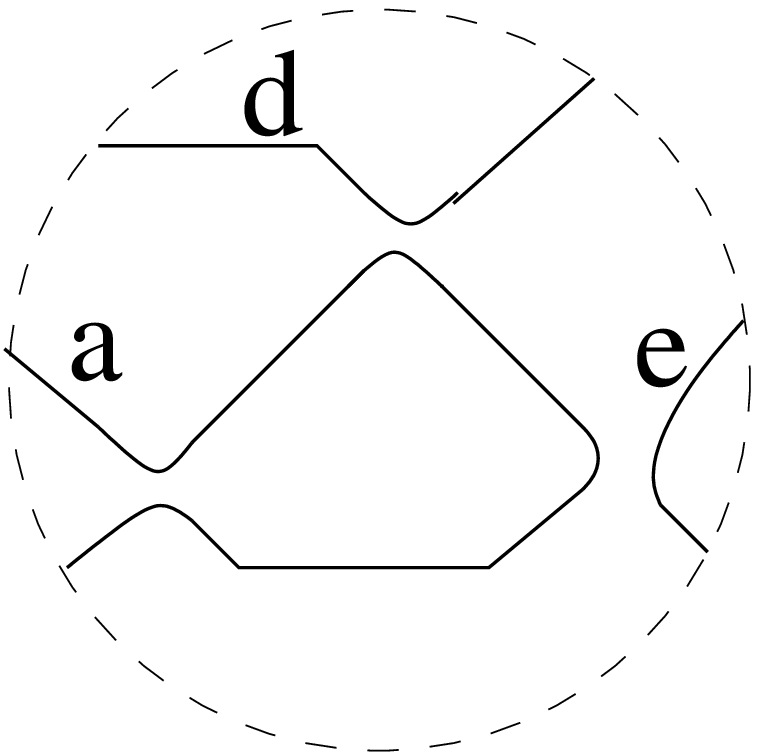,height=1.5cm}},
\end{equation}
where $m(S)$ is the number of negative crossing markers of $S,$
the second sum is as before, and the first sum is over all $f,g$ 
such that \parbox{1.5cm}{\psfig{figure=r3fg.eps,height=1.5cm}} 
$\in C_{i,j+2,s}(D).$
Since the states in the first sum have $-$ and $+$ markers at $v$ and
$w,$ respectively, and the states in the second sum have $+$ and $-$ 
markers at $v$ and $w,$
\begin{equation}\label{f}
f\hat\gamma\rho_{II}(S)=\hat\gamma\rho_{II}(S),
\end{equation}
 by the definition of $f.$
 
Similarly considering $S$ as a state
in $D_{++-}',$ we get
\comment{
$$\rho_{II}'(S)=\parbox{1.5cm}{\psfig{figure=r3-2abc.eps,height=1.5cm}}
+ \sum \parbox{1.5cm}{\psfig{figure=r3-2cfg.eps,height=1.5cm}}.$$
Therefore, by (\ref{gamma}) and Proposition \ref{gamma_prop}(4),
} 
\begin{equation}\label{partialrho'}
\hat\gamma\rho_{II}'(S)=(-1)^{m(S)+1}\sum 
\parbox{1.5cm}{\psfig{figure=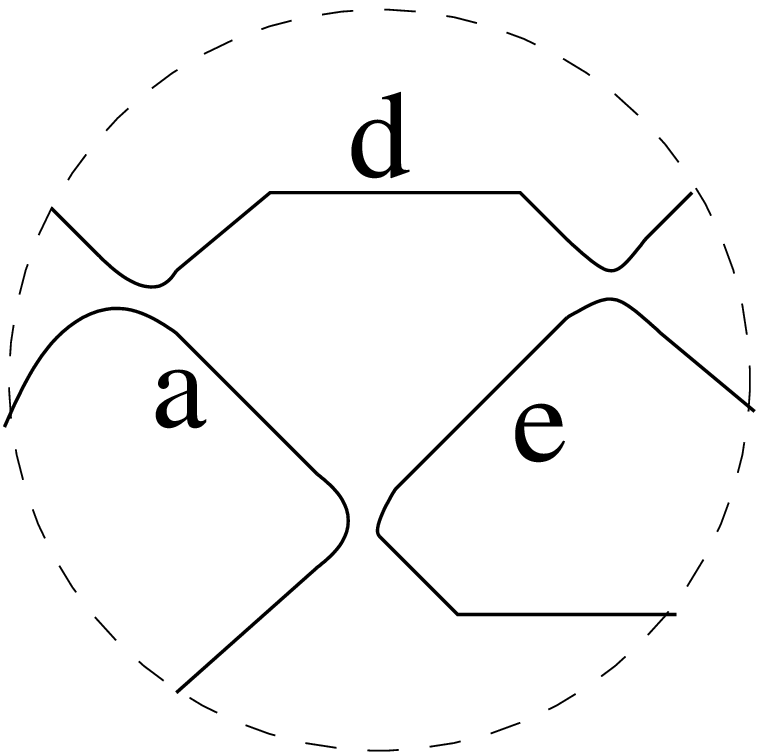,height=1.5cm}}
+ (-1)^{m(S)+1}\sum
\parbox{1.5cm}{\psfig{figure=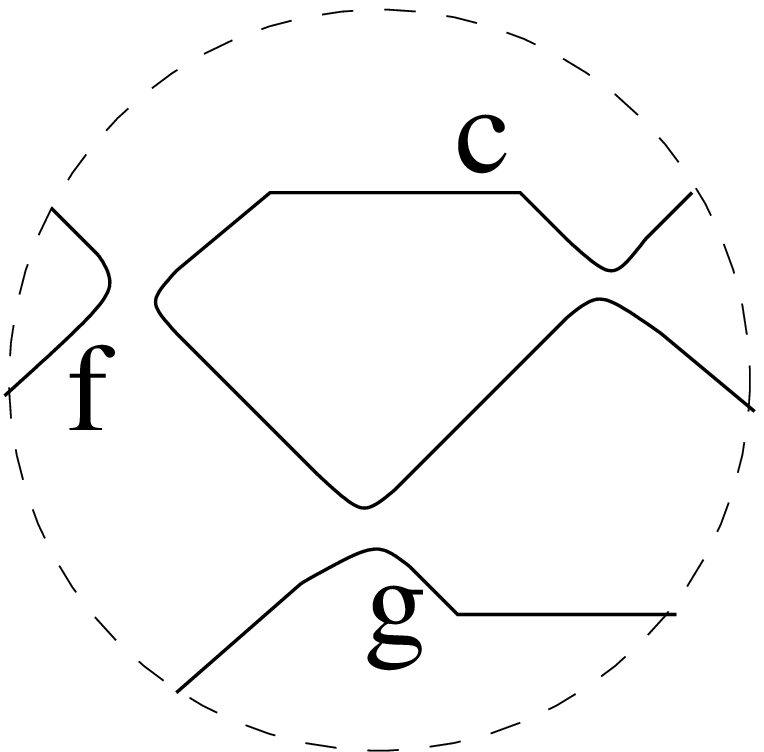,height=1.5cm}},
\end{equation}
by (\ref{gamma}) and Proposition \ref{gamma_prop}(4).

Now (\ref{parrho}),(\ref{f}), and (\ref{partialrho'})
imply (\ref{first_sq}).
\end{proof}

\begin{proposition}\label{IIIchain}
The map $\rho_{III}: C_{ijs}'(D)\to C_{ijs}'(D')$ defined in Proposition
\ref{IIIcommutes} is a chain map.
\end{proposition}

\begin{proof}
If $S\in C_{ijs}'(D)$ then
\begin{equation}\label{IIId}
\rho_{III}d(S)=\bar\beta\rho\be d(S) + \al f \bar\al d(S).
\end{equation}
Since $\beta$ and $\rho$ are chain maps and $\bar \beta d=
\sum_{q\ne p} \hat d_q \bar \beta,$ the first term on the
right equals
\begin{equation}\label{1st_term}
\bar\be \rho\be d(S)=\bar\be d\rho\be(S)=d\bar\be \rho\be(S)-
\hat d_p \bar\be \rho\be(S).
\end{equation}
We can write $S$ as $S_+ +S_-,$ where $S_\pm$ is
composed of states of $D$ with $\pm$ marker at the
crossing $p.$ Since $\bar \alpha$ commutes with $\hat d_q$ for $q\ne
p,$ 
$$\bar \alpha d(S)=\sum_{q\ne p} \hat d_q \bar \alpha(S)+\bar \alpha
\hat d_p(S)=d\bar \alpha(S)+\bar \alpha \hat d_p(S_+).$$
Hence, for the second term on the right side of (\ref{IIId})
we have
\begin{equation}\label{2nd_term}
\al f \bar\al d(S)= d\alpha f \bar \alpha(S)+\alpha f \bar \alpha 
\hat d_p(S_+).
\end{equation}
Therefore, by (\ref{IIId}),(\ref{1st_term}), and (\ref{2nd_term}), 
$$\rho_{III}d(S) = \left(d \bar\beta\rho\be(S_+) - \dh_p\bar \beta \rho
\be(S)\right)- \left(d\al f \bar\al(S) + \al f \bar\al \dh_p(S_+)\right)$$
$$=d \rho_{III}(S) - X,$$
where $$X=\dh_p\bar\beta \rho\be(S) -\al f \bar\al \dh_p(S_+).$$
We are going to show that $X$ vanishes.
Denote $\beta(S_+)$ by $\tilde S_+.$ Since $\hat\gamma =\bar \alpha \hat
d_p\bar \beta,$ 
$$X=\alpha \bar \alpha\dh_p\bar\beta \rho(\tilde S_+) -
\al f \bar\al \dh_p\bar\beta(\tilde S_+)=\alpha \hat\gamma \rho(\tilde
S_+) -\al f \hat\gamma (\tilde S_+).$$
Hence, $X$ vanishes by Lemma \ref{partialcommutes}.
\end{proof}

Propositions \ref{embedding}, \ref{IIIcommutes} and \ref{IIIchain}
imply that the following diagram is exact and commutes:
{\small
$$\begin{array}{c@{}c@{}c@{}c@{}c}
H_{i+1,j-1,s}(D_+) \stackrel{\partial_*}{\to} & H_{i+1,j+1,s}(D_-) 
\stackrel{\al_*}{\to} & H_{ijs}(D) \stackrel{\be_*}{\to} & 
H_{i-1,j-1,s}(D_+) \stackrel{\partial_*}{\to} & H_{i-1,j+1,s}(D_-)\\
 \downarrow \rho_* & \downarrow f_* & \downarrow \rho_{III*} & 
\downarrow \rho_* & \downarrow f_*  \\
H_{i+1,j-1,s}(D_+') \stackrel{\partial_*}{\to} & H_{i+1,j+1,s}(D_-')
\stackrel{\al_*}{\to} & H_{ij}(D') \stackrel{\be_*}{\to} & 
H_{i-1,j-1,s}(D_+') \stackrel{\partial_*}{\to} &
H_{i-1,j+1,s}(D_-').
\end{array} $$}%
By Five Lemma $\rho_{III*}: H_{ijs}(D)\to H_{ijs}(D')$
is an isomorphism. Hence the proof of of Theorem \ref{main_thm}(1) is 
completed.

\np

\Addresses\recd


\begin{thebibliography}

\bibitem[AP]{AP} {\bf M\,M~Asaeda}, {\bf J\,H~Przytycki},
\emph{Khovanov homology: torsion and thickness}, to appear in:
``Advances in Topological Quantum Field Theory (Kananaskis 
Village, Canada 2001)'', (John M. Bryden, editor), 
\arxiv{math.GT/0402402}

\bibitem[APS]{APS}  {\bf M\,M~Asaeda}, {\bf J\,H~Przytycki},
{\bf A\,S Sikora}, 
\emph{A categorification of the skein module of tangles},
\arxiv{math.QA/0410238}

\bibitem[BN]{BN1}
{\bf D Bar-Natan}, \emph{On Khovanov's categorification of the 
Jones polynomial},
\agtref2{2002}{16}{337}{270} \MR{1917056}

\bibitem[H]{Ha}
{\bf A Hatcher}, \emph{Algebraic Topology}, Cambridge University Press (2002)
\url{http://www.math.cornell.edu/~hatcher/AT/ATch3.pdf}

\bibitem[J]{J}
{\bf M Jacobsson},
{\em An invariant of link cobordisms from Khovanov's homology theory},
\agtref4{2004}{53}{1211}{1251}

\bibitem[K1]{K1}
{\bf M Khovanov}, {\em A categorification of the Jones polynomial},
Duke Math. J. 101 (2000) 359--426 \MR{1740682}


\bibitem[K2]{K2}
{\bf M Khovanov}, \emph{Categorifications of the colored Jones polynomial},
\arxiv{math.QA/0302060}

\bibitem[K3]{K3}
{\bf M Khovanov},  \emph{An invariant of tangle cobordisms},
\arxiv{math.QA/0207264}

\bibitem[K4]{K4}
{\bf M Khovanov}, \emph{A functor-valued invariant of tangles},
Algebr. Geom. Topol. \href{http://www.maths.warwick.ac.uk/agt/AGTVol2/agt-2-30.abs.html}{2 (2002) 665--741}
\MR{1928174}

\bibitem[Li]{Li}
{\bf J\,B Listing}, \emph{Vorstudien zur Topologie},  G\"ottinger Studien
(Abtheilung 1) 1 (1847) 811--875

\bibitem[Pr]{Pr}
{\bf J\,H Przytycki}, \emph{Fundamentals of Kauffman bracket skein modules},
Kobe Math. J. 16 (1999) 45--66
\MR{1723531}

\bibitem[PS]{PS} {\bf J\,H Przytycki}, {\bf A\,S Sikora},
{\em On Skein Algebras and $Sl_2(\mathbb C)$-Character Varieties},
Topology 39 (2000) 115--148
\MR{1710996}

\bibitem[V]{V}
{\bf O Viro}, {\em Remarks on definition of Khovanov homology}, 
Fund. Math.
to appear,
\arxiv{math.GT/0202199}

\end{thebibliography}
\end{document}